\documentclass[12pt,oneside,reqno]{amsart}
\usepackage{amssymb}
\usepackage{graphicx}
\author{Alexander Gorodnik}
\title{Uniform distribution of orbits of lattices on spaces of frames}
\address{Department of Mathematics, University of Michigan, Ann Arbor, MI 48109}
\email{gorodnik@umich.edu}
\thanks{This article is a part of author's PhD thesis at Ohio State University done under supervision of Prof.~Bergelson.}

\newtheorem{thm}{Theorem}
\newtheorem{lem}[thm]{Lemma}
\newtheorem{pro}[thm]{Proposition}
\newtheorem{col}[thm]{Corollary}

\newenvironment{example}[1][Example]{\begin{trivlist}
  \item[\hskip \labelsep {\bfseries #1}]}{\end{trivlist}}
\newenvironment{remark}[1][Remark.]{\begin{trivlist}
   \item[\hskip \labelsep {\bfseries #1}]}{\end{trivlist}}
   
\begin{document}
\begin{abstract}
We study distribution of orbits of a lattice $\Gamma\subseteq\hbox{SL}(n,\mathbb{R})$
in the the space $\mathcal{V}_{n,l}$ of $l$-frames in $\mathbb{R}^n$
($1\le l\le n-1$). Examples of dense $\Gamma$-orbits are known from the work
of Dani, Raghavan, and Veech. We show
that dense orbits of $\Gamma$ are uniformly distributed in $\mathcal{V}_{n,l}$ with respect to
an explicitly described measure. We also establish analogous result for lattices in $\hbox{Sp}(n,\mathbb{R})$
that act on the space of isotropic $n$-frames.
\end{abstract}

\maketitle

\section{Introduction}

Let $G=\hbox{SL}(n,\mathbb{R})$ and $\mathcal{V}_{n,l}$ be the space of $l$-frames in
$\mathbb{R}^n$ (i.e. the space of $l$-tuples of linearly independent
vectors in $\mathbb{R}^n$), $1\le l\le n$. The group $G$ acts on this space as follows:
$$
g\cdot (v_1,\ldots, v_l)=(gv_1,\ldots,gv_l),\quad g\in\hbox{SL}(n,\mathbb{R}),\; (v_1,\ldots, v_l)\in \mathcal{V}_{n,l}.
$$
The action is transitive for $l<n$.
Let $\Gamma$ be a lattice in $G$; that is, a discrete subgroup in $G$ such that
the factor space $\Gamma\backslash G$ has finite volume (e.g. $\Gamma=\hbox{SL}(n,\mathbb{Z})$).
The main result of this paper concerns  distribution of $\Gamma$-orbits in $\mathcal{V}_{n,l}$.

When $l=n$, every orbit of $\Gamma$ is discrete. The situation becomes much more
interesting for $l<n$. Let us recall known results:

\begin{thm}{\bf (Dani, Raghavan \cite{dr80})} \label{th_dr}
Let $\Gamma=\hbox{SL}(n,\mathbb{Z})$, and
$v=(v_1,\ldots, v_l)$ be an $l$-frame in $\mathbb{R}^n$, $1\le l\le n-1$.
Then the orbit $\Gamma\cdot v$ is dense in $\mathcal{V}_{n,l}$
iff the space spanned by $\{v_i:i=1,\ldots, l\}$ contains no nonzero rational vectors.
\end{thm}

\begin{thm}{\bf (Veech \cite{ve77})} \label{th_ve}
If $\Gamma$ is a cocompact lattice in $G$, then every orbit of $\Gamma$
in $\mathcal{V}_{n,l}$, $1\le l\le n-1$, is dense.
\end{thm}

Theorems \ref{th_dr} and \ref{th_ve} provide examples of dense $\Gamma$-orbits in $\mathcal{V}_{n,l}$.
Here we show that dense $\Gamma$-orbits
are uniformly distributed with respect to an explicitly described measure on $\mathcal{V}_{n,l}$.
This measure is $\frac{dv}{\hbox{\small \rm Vol}(v)}$, where $dv$ is the Lebesgue measure
on $\prod_{i=1}^{l} \mathbb{R}^n$, and $\hbox{\rm Vol}(v)$ is the $l$-dimensional volume of the frame $v$.

Note that the measure $dv$ is $G$-invariant, and it is unique up to a constant.
However, orbits of $\Gamma$ are equidistributed with respect to the measure $\frac{dv}{\hbox{\small \rm Vol}(v)}$,
which is not $G$-invariant. This phenomenon was already observed by Ledrappier \cite{l}.


Define a norm on $M(n,\mathbb{R})$ by 
\begin{equation}\label{eq_norm}
\|x\|=\left(\hbox{Tr}({}^tx x)\right)^{1/2}=\left(\sum_{i,j} x_{ij}^2\right)^{1/2}\;\;\hbox{for}\;\; x=(x_{ij})\in M(n,\mathbb{R}).
\end{equation}
For $T>0$, $\Omega\subseteq \mathcal{V}_{n,l}$, $v^0\in \mathcal{V}_{n,l}$, put
\begin{equation}\label{eq_N_T}
N_T(\Omega,v^0)=|\{\gamma\in\Gamma: \|\gamma\|<T,\gamma\cdot v^0\in\Omega\}|.
\end{equation}
We determine asymptotic behavior of $N_T(\Omega,v^0)$ as $T\rightarrow\infty$.
This result gives a quantitative strengthening of Theorems \ref{th_dr} and \ref{th_ve},
and it can be interpreted as uniform distribution of dense orbits of $\Gamma$ in $\mathcal{V}_{n,l}$. 

\begin{thm} \label{th_frames00}
Let $\Gamma$ be a lattice in $\hbox{SL}(n,\mathbb{R})$.
Let $v^0\in {\mathcal{V}}_{n,l}$ be an $l$-frame in $\mathbb{R}^n$
such that $\Gamma\cdot v^0$ is dense in ${\mathcal{V}}_{n,l}$.
Let $\Omega$ be a relatively compact Borel subset of ${\mathcal{V}}_{n,l}$ such that $\int_{\partial\Omega}dv=0$.
Then 
\begin{equation} \label{eq_f_main00}
N_T(\Omega,v^0)\sim a_{n,l} \frac{\hbox{\rm Vol}(v^0)^{1-n}}{\bar\mu (\Gamma\backslash G)}\left(\int_\Omega\frac{dv}{\hbox{\small \rm Vol}(v)}\right) T^{(n-1)(n-l)}\quad\hbox{as}\quad T\rightarrow\infty,
\end{equation}
where $a_{n,l}$ is a constant (which is computed in (\ref{eq_anl}) below),
and $\bar \mu$ is a $G$-invariant measure on $\Gamma\backslash G$ (which is defined in (\ref{eq_mubar}) below).
\end{thm}

\begin{remark}
The term $T^{(n-1)(n-l)}$ in (\ref{eq_f_main00}) comes from the asymptotics of the volume of
the set $\{h\in H:\|h\|<T\}$ in the stabilizer $H$ of $v^0$ with respect to the measure on $H$
which is determined by the choice
of the Haar measures on $G$ and $\mathcal{V}_{n,l}=G\cdot v^0$ (see Section \ref{sec_ttt}).
\end{remark}

For $n=2$ and $l=1$, this theorem was proved by Ledrappier \cite{l} for general $\Gamma$ 
and by Nogueira \cite{no} for $\Gamma=\hbox{SL}(2,\mathbb{Z})$ and $\max$-norm using different methods.

Combining Theorems \ref{th_dr} and \ref{th_frames00}, we get:
\begin{col} \label{th_frames}
Let $\Gamma=\hbox{SL}(n,\mathbb{Z})$.
Let $v^0=(v_1^0,\ldots, v_{l}^0)\in {\mathcal{V}}_{n,l}$ be an $l$-frame in $\mathbb{R}^n$
such that the space $\left<v^0_1,\ldots, v^0_{l}\right>$ contains no nonzero rational vectors.
Let $\Omega$ be a relatively compact Borel subset of ${\mathcal{V}}_{n,l}$ such that $\int_{\partial\Omega}dv=0$.
Then 
\begin{equation} \label{eq_f_main}
N_T(\Omega,v^0)\sim b_{n,l} \hbox{\rm Vol}(v^0)^{1-n}\left(\int_\Omega \frac{dv}{\hbox{\small \rm Vol}(v)}\right) T^{(n-1)(n-l)}\quad\hbox{as}\quad T\rightarrow\infty,
\end{equation}
where $b_{n,l}$ is a constant (which is computed in (\ref{eq_bnl}) below).
\end{col}

\begin{example}
Figure \ref{pic1} shows a part of the the orbit $\hbox{SL}(2,\mathbb{Z})v^0$
for $v^0={}^t (\sqrt{2},\sqrt{3})$. By the result of Ledrappier,
this orbit is uniformly distributed in $\mathbb{R}^2$ with respect to the measure
$\frac{dxdy}{\sqrt{x^2+y^2}}$.
\end{example}

\begin{figure}[h] \label{pic1}
\begin{center}
\includegraphics[width=7cm,height=7cm]{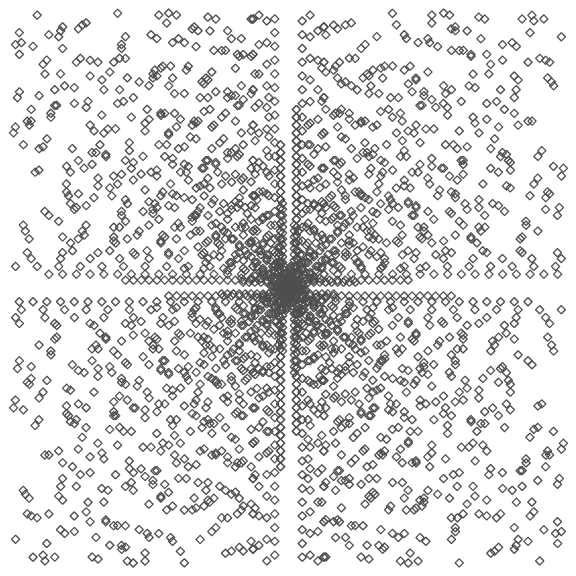}
\end{center}
\caption{}
\end{figure}

Dani and Raghavan also considered orbits of frames under $\hbox{Sp}(n,\mathbb{Z})$.
Denote
$$
J=\left(
\begin{tabular}{c|c}
$0$ & $E$\\
\hline
$-E$ & $0$
\end{tabular}
\right),
$$
where $E$ is the identity $n\times n$ matrix. The symplectic form 
$(x,y)\mapsto {}^txJy$ will be denoted by $J$ too.
Let
$$
G=\hbox{Sp}(n,\mathbb{R})=\{g\in\hbox{SL}(2n,\mathbb{R}): {}^tgJg=J\}
$$
and $\Gamma=\hbox{Sp}(n,\mathbb{Z})$.
A frame $(v_1,\ldots, v_s)$ is called {\it isotropic} if the symplectic form $J$ is
$0$ on the space spanned by $\{v_i:i=1,\ldots,s\}$.

\begin{thm} {\bf (Dani, Raghavan \cite{dr80})} \label{th_dr2}
Let $v=(v_1,\ldots, v_{n})$ be an isotropic frame in $\mathbb{R}^{2n}$.
Then $\Gamma\cdot v$ is dense in the space
of isotropic $n$-frames iff the space spanned by $\{v_i:i=1,\ldots, n\}$
contains no nonzero rational vectors.
\end{thm}

A result similar to Theorem \ref{th_frames00} holds in this case too.
Denote by ${\mathcal{W}}_n$ the space of $2n$-dimensional $n$-frames that are isotropic with respect to the standard
symplectic form $J$. Note that ${\mathcal{W}}_n$ is an open subset of an algebraic set in $\prod_{i=1}^n \mathbb{R}^{2n}$.
Since by Witt's theorem $\hbox{Sp}(n,\mathbb{R})$ acts transitively on ${\mathcal{W}}_n$,
${\mathcal{W}}_n$ is a submanifold of $\prod_{i=1}^n \mathbb{R}^{2n}$.

We improve Theorem \ref{th_dr2} by showing that dense orbits of $\Gamma$ are uniformly distributed:

\begin{thm} \label{th_frames0}
Let $\Gamma$ be a lattice in $\hbox{Sp}(n,\mathbb{R})$,
and $v^0\in {\mathcal{W}}_n$ be such that $\Gamma \cdot v^0$ is dense in ${\mathcal{W}}_n$.
Let $\Omega$ be a relatively compact Borel subset of ${\mathcal{W}}_n$ such that the boundary of $\Omega$
has measure $0$ in the Lebesgue measure class.
Then 
$$
N_T(\Omega,v^0)\sim \lambda_{v^0}(\Omega)T^{n(n+1)/2}\quad\hbox{as}\quad T\rightarrow\infty
$$
for some measure $\lambda_{v^0}$ on ${\mathcal{W}}_n$ in the Lebesgue measure class,
which can be explicitly computed.
\end{thm}

Note that the measure $\lambda_{v^0}$ is not $\hbox{Sp}(n,\mathbb{R})$-invariant.

In the next section we show how to derive asymptotic distribution for counting functions
similar to $N_T(\Omega,v^0)$ from uniform distribution of orbits of  subgroups of $G$ in the space
$\Gamma\backslash G$. In section \ref{sec_sl}, we consider the case $G=\hbox{SL}(n,\mathbb{R})$.
First, for convenience of the reader, we sketch an easy proof of Theorems \ref{th_dr} and \ref{th_ve}
based on topological rigidity of unipotent flows, which was established by Ratner \cite{r91b}.
Then we introduce a decomposition of $G$ based on the Iwasawa decomposition, and obtain
results on volume of balls in the subgroup $B^o_l$, which is defined below.
This allows us to use results from Section \ref{sec_ttt} to prove
Theorem \ref{th_frames00} and Corollary \ref{th_frames} modulo ergodic theorem 
along balls in the group $B^o_l$ (Theorem \ref{th_ergodic}).
In Section \ref{sec_erg}, we prove the ergodic theorem for $B^o_l$.
Note that for $l=n-1$ it is a special case of the result of Shah \cite{sh94}.
The proof of the ergodic theorem is similar to the proof of Theorem 2 from \cite{go}.
Finally, in Section \ref{sec_sp} we consider the case $G=\hbox{Sp}(n,\mathbb{R})$, and
prove Theorem \ref{th_frames0}. The method of the proof is similar to one used
for Theorem \ref{th_frames00}: we use Iwasawa decomposition for $\hbox{Sp}(n,\mathbb{R})$
and uniform distribution of large unipotent subgroups due to Shah \cite{sh94}.
In the Appendix, we prove some technical volume estimates and Corollary \ref{th_frames}.

\begin{remark}
In the definition of $N_T(\Omega,v^0)$, we used the norm (\ref{eq_norm}).
The fact that this norm is invariant under the orthogonal group made our calculations easier.
However, one can prove similar results for every norm on $\hbox{M}(n,\mathbb{R})$
with possibly different limit measure in the Lebesgue measure class.
\end{remark}

\noindent {\sc Acknowledgment:} I would like to thank H.~Furstenberg for fruitful discussions and
Barak Weiss for some comments and for pointing out an error in the preliminary version of this paper.
I am also very grateful to V. Bergelson for his encouragement and for many interesting discussions.

\section{Some limit theorems} \label{sec_ttt}

In this section we establish asymptotics of some counting functions.

Let $G\subseteq \hbox{SL}(n,\mathbb{R})$ be a Lie group, $\Gamma$ a
lattice in $G$, and $H$ a Lie subgroup of $G$. Denote by $\varrho$ a right Haar measure on $H$.
Let $\mu$ be a Haar measure on $G$, and $\bar\mu$ be the measure on $\Gamma\backslash G$
such that
$$
\int_G fd\mu=\int_{\Gamma\backslash G} \left(\sum_{\gamma\in\Gamma} f(\gamma g)\right)d\bar\mu(g),\quad f\in
C_c(G).
$$
Throughout this section, we assume that
for some $M>0$ and every $c>0$, $m\in\mathbb{R}$,
\begin{equation}\label{eq_H1}
\lim_{T\rightarrow\infty}\frac{\varrho(H_{cT+m})}{\varrho(H_T)}=c^M.
\end{equation}
where $H_T=\{h\in H:\|h\|<T\}$, and for every $x\in \Gamma\backslash G$ such that $xH$ is dense in 
$\Gamma\backslash G$ and $\tilde{f}\in C_c(\Gamma\backslash G)$,
\begin{equation}\label{eq_H2}
\frac{1}{\varrho(H_T)}\int_{H_T} \tilde{f}(x h^{-1})d\varrho(h)\rightarrow \frac{1}{\bar\mu(\Gamma\backslash G)}\int_{\Gamma\backslash G} \tilde{f}d\mu\quad\hbox{as}\quad T\rightarrow\infty.
\end{equation}
 
First, we prove an elementary lemma:

\begin{lem} \label{hyp_l2}
Let $(V,\|\cdot\|)$ be a normed vector space, $G$ be a topological group, and 
$\rho: G\rightarrow \mathcal{B}(V)^*$ ($=$ the space of bounded invertible linear operators on $V$) be a continuous map (w.r.t. norm topology).
Then for any $g_0\in G$ and $k>1$, there exists a neighborhood ${\mathcal{O}}_{g_0}$ of $g_0$ in $G$ such that
for any $g\in {\mathcal{O}}_{g_0}$ and $v\in V$, 
$$
k^{-1}\|\rho(g_0)v\|\le \|\rho(g)v\|\le k\|\rho(g_0)v\|.
$$
\end{lem}

\begin{proof} 
Take 
$$
{\mathcal{O}}_{g_0}=\{g:\|\rho(g_0)\rho(g)^{-1}\|<k,\|\rho(g)\rho(g_0)^{-1}\|<k\}.
$$
Then for $g\in {\mathcal{O}}_{g_0}$, 
$$
\|\rho(g_0)v\|\le \|\rho(g_0)\rho(g)^{-1}\|\cdot\|\rho(g)v\|<k\|\rho(g)v\|.
$$
Similarly, for $g\in {\mathcal{O}}_{g_0}$, 
$$
\|\rho(g)v\|\le\|\rho(g)\rho(g_0)^{-1}\|\cdot\|\rho(g_0)v\|<k\|\rho(g_0)v\|.
$$
\end{proof}

For $g\in\hbox{SL}(n,\mathbb{R})$, denote by $\hat g:x\mapsto gxg^{-1}$ the inner automorphism
corresponding to $g$. For a subgroup $L$ of $G$, denote by $N_G(L)$ the normalizer of $L$ in $G$.
For $g\in N_G(H)$, define
\begin{equation}\label{eq_DH}
\Delta_H(g)=\left|\det \left(\hbox{\rm Ad}(g)|_{{\hbox{\rm\small Lie}(H)}}\right)\right|
\end{equation}
where $\hbox{\rm Ad}(g)$ is the adjoint transformation of the Lie algebra of $G$.

\begin{pro} \label{lem_fl}
Let $x_0\in G$ be such that $\Gamma x_0 H$ is dense in $\Gamma\backslash G$.
Let
\begin{eqnarray*}
g\mapsto c_g:&& G\rightarrow N_G(H),\\
g\mapsto d_g:&& G\rightarrow H,\\
g\mapsto e_g:&& G\rightarrow \mathbb{R},\\
g\mapsto m_g:&& G\rightarrow \mathbb{R}_+,\\
\end{eqnarray*}
be continuous maps that factor through $G/\hat x_0(H)$.
Then for every $f\in C_c(G)$,
\begin{eqnarray} \label{eq_fl}
\mathop{\lim}_{T\rightarrow\infty}\frac{1}{\varrho(H_T)}\sum_{\gamma\in\Gamma}\int_{m_\gamma^2\|\hat c_\gamma (hd_\gamma)\|^2+e_\gamma<T^2} f(\gamma x_0 h^{-1}x_0^{-1})d\varrho(h)\\
=\frac{1}{\bar \mu(\Gamma\backslash G)}\int_{G} \frac{f(g)}{m_g^M\cdot
\Delta_H(c_g)}d\mu(g).\nonumber
\end{eqnarray}
\end{pro}

\begin{proof}
We shall assume without loss of generality that  $f\ge 0$.

There exist real $M_1$ and $M_2$ such that $M_1\le e_g\le M_2$ for $g\in\hbox{supp}(f)\hat x_0(H)$.
Then for $g\in\hbox{supp}(f)\hat x_0(H)$, 
\begin{equation}\label{eq_lim_est}
\{h:m_g^2\|\hat c_g(hd_g)\|^2+e_g<T^2\}\subseteq \{h:m_g\|\hat c_g(hd_g)\|<(T^2-M_1)^{1/2}\}.
\end{equation}

Denote $\tilde{f}(\Gamma g)=\sum_{\gamma\in \Gamma}f(\gamma g)$.
Note that $\tilde{f}\in C_c(\Gamma\backslash G)$. 

Let $r>1$ and $\varepsilon>0$. By Lemma \ref{hyp_l2}, for any $g_0\in G$ 
there exists a neighborhood ${\mathcal{O}}_{g_0}$ of $g_0$ such that
\begin{eqnarray} 
\label{eq_lim_est1} r^{-1}\|\hat c_{g_0}(v)\|<\|\hat c_g(v)\|<r\|\hat c_{g_0}(v)\|,\\
\label{eq_lim_est10} r^{-1}\|v\hat c_{g_0}(d_{g_0})\|<\|v\hat c_{g_0}(d_g)\|<r\|v\hat c_{g_0}(d_{g_0})\|
\end{eqnarray}
for all $g\in {\mathcal{O}}_{g_0}\hat x_0(H)$ and $v\in M(n,\mathbb{R})$.
Moreover, ${\mathcal{O}}_{g_0}$ can be taken such that 
\begin{eqnarray*}
\left| \frac{1}{m_g^M\cdot\Delta_H(c_g)}-\frac{1}{m_{g_0}^M\cdot\Delta_H(c_{g_0})}\right|<\varepsilon,
\end{eqnarray*}
and
\begin{eqnarray*}
r^{-1}m_{g_0}\le m_g\le rm_{g_0}
\end{eqnarray*}
for all $g\in {\mathcal{O}}_{g_0}\hat x_0(H)$.

Note that for every $v\in N_G(H)$, $\Gamma x_0vH$ is dense in $\Gamma\backslash G$. Therefore, by (\ref{eq_H2}),
for every $u\in G$ and $v\in N_G(H)$,
\begin{eqnarray}
\mathop{\lim}_{T\rightarrow\infty}\frac{1}{\varrho(H_T)}\int_{H_T}\tilde{f}(\Gamma x_0 v h^{-1}u)d\varrho (h)
&=&\frac{1}{\bar\mu(\Gamma\backslash G)}\int_{\Gamma\backslash G} \tilde{f}(yu)d\bar\mu(y)\nonumber\\
&=&\frac{1}{\bar\mu(\Gamma\backslash G)}\int_{\Gamma\backslash G} \tilde{f}d\bar\mu.\label{eq_H3}
\end{eqnarray}

To prove (\ref{eq_fl}), we first suppose that $\hbox{supp}(f)\subseteq{\mathcal{O}}_{g_0}$
for some $g_0\in G$. Put $c_0=c_{g_0}$, $d_0=d_{g_0}$, $m_0=m_{g_0}$.
Then using (\ref{eq_lim_est}), (\ref{eq_lim_est1}),  and (\ref{eq_lim_est10}), we get
\begin{eqnarray*}
&&\sum_{\gamma\in\Gamma}\int_{m_\gamma^2\|\hat c_\gamma (hd_\gamma)\|^2+e_\gamma<T^2} f(\gamma x_0h^{-1}x_0^{-1})d\varrho(h)\\
&\le&\sum_{\gamma\in\Gamma}\int_{m_\gamma\|\hat c_\gamma(hd_\gamma)\|<(T^2-M_1)^{1/2}} f(\gamma x_0h^{-1}x_0^{-1})d\varrho(h)\\
&\le&\sum_{\gamma\in\Gamma}\int_{r^{-2}m_0\|\hat c_0(h) \hat c_0(d_\gamma)\|<(T^2-M_1)^{1/2}} f(\gamma x_0h^{-1}x_0^{-1})d\varrho(h)\\
&\le&\sum_{\gamma\in\Gamma}\int_{r^{-3}m_0\|\hat c_0(h) \hat c_0(d_0)\|<(T^2-M_1)^{1/2}} f(\gamma x_0h^{-1}x_0^{-1})d\varrho(h)\\
&=&\int_{\|h\|<r^3m_0^{-1}(T^2-M_1)^{1/2}} \tilde{f}(\Gamma x_0c_0^{-1}\hat c_0(d_0) h^{-1}c_0 x_0^{-1})\frac{d\varrho(h)}{\Delta_H(c_0)}.
\end{eqnarray*}
Thus, by (\ref{eq_H3}) and (\ref{eq_H1}),
\begin{eqnarray*}
&&\limsup_{T\rightarrow\infty}\frac{1}{\varrho(H_T)}\sum_{\gamma\in\Gamma}\int_{m_\gamma^2\|\hat c_\gamma (hd_\gamma)\|^2+e_\gamma<T^2} f(\gamma x_0h^{-1}x_0^{-1})d\varrho(h)\\
&\le& \limsup_{T\rightarrow\infty} \frac{\varrho(H_{r^3m_0^{-1}(T^2-M_1)^{1/2}})}{\varrho(H_T)} \cdot \frac{1}{\bar \mu(\Gamma\backslash G)\Delta_H(c_0)}\int_{\Gamma\backslash G} \tilde{f}d\bar\mu\\
&=&\frac{r^{3M}}{m_0^M\cdot\bar \mu(\Gamma\backslash G)\cdot\Delta_H(c_0)}\int_{\Gamma\backslash G}
\tilde{f}d\bar\mu\\
&=&\frac{r^{3M}}{m_0^M\cdot\bar \mu(\Gamma\backslash G)\cdot\Delta_H(c_0)}\int_{G} f d\mu.
\end{eqnarray*}

Now let $f$ be arbitrary. There exists a finite cover
$\hbox{supp}(f)\subseteq \cup_i {\mathcal{O}}_{g_i}$. Let $c_i=c_{g_i}$ and $m_i=m_{g_i}$.
Let $\alpha_i\in C_c(G)$ be a partition of unity for $\{{\mathcal{O}}_{g_i}\}$
such that $\sum_i\alpha_i=1$ on $\hbox{supp}(f)$. Then
\begin{eqnarray*}
\limsup_{T\rightarrow\infty}\frac{1}{\varrho(H_T)}\sum_{\gamma\in\Gamma}\int_{m_\gamma^2\|\hat c_\gamma (hd_\gamma)\|^2+e_\gamma<T^2} f(\gamma x_0h^{-1}x_0^{-1})d\varrho(h)\\
=\limsup_{T\rightarrow\infty}\frac{1}{\varrho(H_T)}\sum_{\gamma\in\Gamma}\int_{m_\gamma^2\|\hat c_\gamma (hd_\gamma)\|^2+e_\gamma<T^2} \Big\{ \sum_i f(\gamma x_0h^{-1}x_0^{-1}) \\
\cdot\alpha_i(\gamma x_0h^{-1}x_0^{-1})\Big\} d\varrho(h)
\le \frac{r^{3M}}{\bar \mu(\Gamma\backslash G)} \sum_i\int_G \frac{f(g)}{m_i^M\Delta_H(c_i)}\alpha_i(g)d\mu(g)\\
\le \frac{r^{3M}}{\bar \mu(\Gamma\backslash G)} \sum_i\int_G \left(\frac{f(g)}{m_g^M\cdot\Delta_H(c_g)}+\varepsilon f(g)\right)\alpha_i(g)d\mu(g)\\
=\frac{r^{3M}}{\bar \mu(\Gamma\backslash G)}\int_G \frac{f(g)}{m_g^M\cdot\Delta_H(c_g)}d\mu(g)+
\frac{r^{3M}\varepsilon}{\bar \mu(\Gamma\backslash G)}\int_G fd\mu
\end{eqnarray*}
for every $r>1$ and $\varepsilon>0$. Therefore,
\begin{eqnarray}
\limsup_{T\rightarrow\infty} \frac{1}{\varrho(H_T)}\sum_{\gamma\in\Gamma}\int_{m_\gamma^2\|\hat c_\gamma (hd_\gamma)\|^2+e_\gamma<T^2} f(\gamma x_0hx_0^{-1})d\varrho(h)\nonumber\\
\le\frac{1}{\bar \mu(\Gamma\backslash G)}\int_G \frac{f(g)}{m_g^M\cdot\Delta_H(c_g)}d\mu(g).\label{eq_limfl}
\end{eqnarray}

Similarly, one can prove the lower bound for (\ref{eq_fl}).
\end{proof}

\begin{pro} \label{lem_fl1}
Let $f$ be the characteristic function of a relatively compact Borel subset $Z\subseteq G$
such that $\mu (\partial Z)=0$. Let $x_0\in G$ be such that $\Gamma x_0 H$ is dense in
$\Gamma\backslash G$. Then (\ref{eq_fl}) holds for $f$.
\end{pro}

\begin{proof}
Denote by $\hbox{int}(Z)$ and $\overline{Z}$ the interior and the closure of $Z$ respectively.

Let $W_0$ be an open relatively compact subset such that $\overline{Z}\subseteq W_0$.
There exists $C>0$ such that $\left(m_g^M\cdot\Delta_H(c_g)\right)^{-1}\le C$ for $g\in W_0$.

Let $\varepsilon>0$. There exist a compact subset $V\subseteq \hbox{int}(Z)$
and an open subset $W$ such that $\overline{Z}\subseteq W\subseteq W_0$ and $\mu (W-V)<\varepsilon$.
Take functions $f_1,f_2\in C_c(G)$ such that $0\le f_i\le 1$,
$f_1=1$ on $V$, $f_1=0$ outside $\hbox{int}(Z)$, $f_2=1$ on $\overline{Z}$, and $f_2=0$ outside $W$.
Then $f_1\le f\le f_2$. By Proposition \ref{lem_fl} applied to $f_2$,
\begin{eqnarray*}
&&\limsup_{T\rightarrow\infty}\frac{1}{\varrho(H_T)}\sum_{\gamma\in\Gamma}\int_{m_\gamma^2\|\hat c_\gamma (hd_\gamma)\|^2+e_\gamma<T^2} f(\gamma x_0h^{-1}x_0^{-1})d\varrho(h)\\
&\le& \frac{1}{\bar \mu(\Gamma\backslash G)}\int_{G} \frac{f_2(g)}{m_g^M\cdot\Delta_H(c_g)}d\mu(g)\\
&\le& \frac{1}{\bar \mu(\Gamma\backslash G)}\left(\int_{G} \frac{f(g)}{m_g^M\cdot\Delta_H(c_g)}d\mu(g)
+\int_{G} \frac{f_2(g)-f_1(g)}{m_g^M\cdot\Delta_H(c_g)}d\mu(g)\right)\\
&\le& \frac{1}{\bar \mu(\Gamma\backslash G)}\int_{G} \frac{f(g)}{m_g^M\cdot\Delta_H(c_g)}d\mu(g)
+\frac{C\mu(W-V)}{\bar \mu(\Gamma\backslash G)}\\
&\le& \frac{1}{\bar \mu(\Gamma\backslash G)}\int_{G} \frac{f(g)}{m_g^M\cdot \Delta_H(c_g)}d\mu(g)+\frac{C\varepsilon}{\bar \mu(\Gamma\backslash G)}
\end{eqnarray*}
for every $\varepsilon>0$. This shows (\ref{eq_limfl}).
The dual inequality for $\liminf$ can be proved similarly.
\end{proof}

Suppose that for a closed subset $Y$ of $G$,
the product map $Y\times \hat x_0 (H)\rightarrow G$ be a homeomorphism. 
For $g\in G$, define $y_g\in Y$ and $h_g\in H$ such that $g=y_g\hat x_0(h_g)$.
The map
$$
\alpha:y\mapsto y\cdot\hat x_0(H):Y \rightarrow G/\hat x_0(H).
$$
is a homeomorphism too.
Let $\nu_1$ be a measure on $Y$ such that
\begin{equation} \label{eq_nu}
\int_G fd\mu=\int_{Y}\int_{H} f(y\hat x_0(h))d\varrho(h)d\nu_1(y),\;\; f\in C_c(G).
\end{equation}
Note that such a measure exists because $\mu$ and $\varrho$ are right invariant.
Let $\nu$ be the measure on $G/\hat x_0(H)$ which is the image of $\nu_1$ under $\alpha$, i.e.
\begin{equation} \label{eq_nu1}
\int_Y f(\alpha(y))d\nu_1(y)=\int_{G/\hat x_0(H)}fd\nu,\;\; f\in C_c(G/\hat x_0(H)).
\end{equation}
Note that the measure $\nu$ depends on the choice of the section $Y$.

\begin{pro} \label{pro_assym}
Use notations as in Proposition \ref{lem_fl}.
Let $\Omega$ be relatively compact Borel subset of $G/\hat x_0(H)$ such that $\nu(\partial\Omega)=0$.
Let 
$$
N_T(\Omega)=|\{\gamma\in\Gamma:m_\gamma^2\|\hat c_\gamma (h_\gamma d_\gamma)\|^2+e_\gamma<T^2,\gamma \cdot\hat x_0(H)\in\Omega)\}|.
$$
Then
$$
\frac{N_T(\Omega)}{\varrho(H_T)}\sim \frac{1}{\bar\mu(\Gamma\backslash G)}\int_{\alpha^{-1}(\Omega)}\frac{1}{m_y^M\cdot\Delta_H(c_y)}d\nu_1(y)
$$
as $T\rightarrow\infty$.
\end{pro}

\begin{proof}
Let
$$
\mathcal{O}_\varepsilon=\{h\in H:\|h-1\|<\varepsilon,\|h^{-1}-1\|<\varepsilon\}
$$
for $\varepsilon>0$.

Let $\phi$ be the characteristic function of $\alpha^{-1}(\Omega)\subseteq Y$.
Take $\psi_\varepsilon$ to be the characteristic function of $\mathcal{O}_\varepsilon$
normalized so that $\int_H \psi_\varepsilon d\varrho=1$.
Let 
\begin{equation} \label{eq_f_fl_def}
f_\varepsilon(g)=\phi(y_g)\psi_\varepsilon (h_g)\quad\hbox{for}\; g\in G.
\end{equation}
Note that $f$ satisfies conditions of Proposition \ref{lem_fl1}, but before applying
this proposition, we need a lemma.

\begin{lem} \label{lem_flinq}
For every $r>1$, there exists $\varepsilon>0$ such that
\begin{equation}\label{eq_NTl}
N_{r^{-1}T}(\Omega)\le\sum_{\gamma\in\Gamma}\int_{m_\gamma^2\|\hat c_\gamma (hd_\gamma)\|^2+e_\gamma<T^2} f_\varepsilon(\gamma x_0h^{-1}x_0^{-1})d\varrho(h)\le N_{rT}(\Omega).
\end{equation}
\end{lem}

\begin{proof}
Note that $f_\varepsilon(\gamma x_0h^{-1}x_0^{-1})=0$ for all $h\in H$ unless 
\begin{equation} \label{eq_inom}
y_\gamma\in \alpha^{-1}(\Omega),
\end{equation}
and if the above condition holds,
\begin{eqnarray}
&&\int_{m_\gamma^2\|\hat c_\gamma (hd_\gamma)\|^2+e_\gamma<T^2} f_\varepsilon(\gamma x_0h^{-1}x_0^{-1})d\varrho(h)\nonumber\\
&=&\int_{m_\gamma^2\|\hat c_\gamma (hd_\gamma)\|^2+e_\gamma<T^2} \psi_\varepsilon (h_\gamma h^{-1})d\varrho(h)\nonumber\\
&=&\int_{m_\gamma^2\|\hat c_\gamma (hh_\gamma d_\gamma)\|^2+e_\gamma<T^2} \psi_\varepsilon (h^{-1})d\varrho(h)\nonumber\\
&=&\int_{m_\gamma^2\|\hat c_\gamma (h) \hat c_\gamma(h_\gamma d_\gamma)\|^2+e_\gamma<T^2} \psi_\varepsilon (h)d\varrho(h). \label{inom1}
\end{eqnarray}
Let
$$
I_\gamma\stackrel{def}{=}\int_{m_\gamma^2\|\hat c_\gamma (h) \hat c_\gamma(h_\gamma d_\gamma)\|^2+e_\gamma<T^2} \psi_\varepsilon (h)d\varrho(h).
$$

For $\gamma$ as in (\ref{eq_inom}), there exists $C>0$ such that
$$
\|\hat c_\gamma (v)\|\le C\|v\|\;\;\hbox{and}\;\;\|\hat c_\gamma^{-1}(v)\|\le C\|v\|\;\;\hbox{for all}\;\; v\in \hbox{M}(n,\mathbb{R}).
$$
It follows that for every $\varepsilon>0$,
$$
\mathcal{O}_{\varepsilon/C}\subseteq \hat c_\gamma\left(\mathcal{O}_{\varepsilon}\right)\subseteq\mathcal{O}_{C\varepsilon}.
$$
Therefore, by Lemma \ref{hyp_l2}, there exists $\varepsilon>0$ such that
\begin{equation}\label{eq_cgn}
r^{-1}\|v\|\le \|\hat c_\gamma(h)v\|\le r\|v\|
\end{equation}
for every $v\in \hbox{M}(n,\mathbb{R})$, $h\in \mathcal{O}_\varepsilon$,
and $\gamma\in\Gamma$ such that (\ref{eq_inom}) holds.

Let $\gamma\in\Gamma$ be such that
$$
m_\gamma^2\|\hat c_\gamma(h_\gamma d_\gamma)\|^2+e_\gamma<r^{-2}T^2,
$$
and (\ref{eq_inom}) holds. Then by (\ref{eq_cgn}),
$$
m_\gamma^2\|\hat c_\gamma(h)\hat c_\gamma(h_\gamma d_\gamma)\|^2+e_\gamma\le r^2m_\gamma^2\|\hat c_\gamma(h_\gamma d_\gamma)\|^2+e_\gamma<T^2
$$
for $h\in\mathcal{O}_\varepsilon$.
It follows that the $I_\gamma=1$. This proves the first inequality in (\ref{eq_NTl}).

Note that $I_\gamma\le 1$. Let $\gamma\in\Gamma$ be such that $I_\gamma\ne 0$.
Then (\ref{eq_inom}) holds, and for some $h\in\mathcal{O}_\varepsilon$,
$$
m_\gamma^2\|\hat c_\gamma(h)\hat c_\gamma(h_\gamma d_\gamma)\|^2+e_\gamma<T^2.
$$
Using (\ref{eq_cgn}), we deduce that
$$
m_\gamma^2\|\hat c_\gamma(h_\gamma d_\gamma)\|^2+e_\gamma<r^{2}T^2.
$$
This proves the second inequality in (\ref{eq_NTl}).
\end{proof}

Now we can use Proposition \ref{lem_fl1} to find asymptotics for $N_T(\Omega)$.
By Lemma \ref{lem_flinq}, for every $r>1$ there exists $\varepsilon>0$ such that 
\begin{eqnarray*}
&&\limsup_{T\rightarrow\infty}\frac{N_T(\Omega)}{\varrho(H_T)}\\
&\le& \limsup_{T\rightarrow\infty} \frac{1}{\varrho(H_T)}\sum_{\gamma\in\Gamma}\int_{m_\gamma^2\|\hat c_\gamma(hd_\gamma)\|^2+e_\gamma<r^2T^2} f_\varepsilon(\gamma x_0h^{-1}x_0^{-1})d\varrho(h).
\end{eqnarray*}
Therefore, by Proposition \ref{lem_fl1}, (\ref{eq_H1}), (\ref{eq_nu}), and (\ref{eq_nu1}),
\begin{eqnarray*}
&&\limsup_{T\rightarrow\infty}\frac{N_T(\Omega)}{\varrho(H_T)}\\
&\le&\left(\limsup_{T\rightarrow\infty} \frac{\varrho(H_{rT})}{\varrho(H_T)}\right)
\frac{1}{\bar\mu(\Gamma\backslash G)}\int_G\frac{f_\varepsilon(g)}{m_g^M\cdot \Delta_H(c_g)}d\mu(g)\\
&=&\frac{r^M}{\bar\mu(\Gamma\backslash G)}\int_{Y\times H}\frac{f_\varepsilon(y\hat x_0(h))}{m_y^M\cdot\Delta_H(c_y)}d\nu_1(y)d\varrho(h)\\
&=&\frac{r^M}{\bar\mu(\Gamma\backslash G)}\int_Y\frac{\phi(y)}{m_y^M\cdot\Delta_H(c_y)}d\nu_1(y)\int_H \psi_\varepsilon(h)d\varrho(h)\\
&=&\frac{r^M}{\bar\mu(\Gamma\backslash G)}\int_{\alpha^{-1}(\Omega)}\frac{d\nu_1(y)}{m_y^M\cdot\Delta_H(c_y)}.
\end{eqnarray*}
Taking $r\rightarrow 1+$, we get
$$
\limsup_{T\rightarrow\infty}\frac{N_T(\Omega)}{\varrho(H_T)}\le
\frac{1}{\bar\mu(\Gamma\backslash G)}\int_{\alpha^{-1}(\Omega)}\frac{d\nu_1(y)}{m_y^M\cdot\Delta_H(c_y)}.
$$
Similarly, one can prove that
$$
\liminf_{T\rightarrow\infty}\frac{N_T(\Omega)}{\varrho(H_T)}\ge
\frac{1}{\bar\mu(\Gamma\backslash G)}\int_{\alpha^{-1}(\Omega)}\frac{d\nu_1(y)}{m_y^M\cdot\Delta_H(c_y)}.
$$
This proves the Proposition.
\end{proof}

\section{Uniform distribution for a lattice in $\hbox{SL}(n,\mathbb{R})$}\label{sec_sl}

\subsection{Density of orbits}\label{sec_dense}

\hspace{1cm}\vspace{0.2cm}

In this section we derive Theorems \ref{th_dr} and \ref{th_ve} from the following result
on topological rigidity of unipotent flow, which was proved by M.~Ratner:

\begin{thm}\label{th_rat}
{\bf (Ratner \cite{r91b})}
Let $G$ be a connected Lie group, $\Gamma$ be a lattice in $G$, and
$U$ be a subgroup of $G$ generated by $\textrm{Ad}$-unipotent $1$-parameter subgroups.
Then for every $x\in\Gamma\backslash G$, $\overline{xU}=xH$, where $H$ is
a closed connected subgroup of $G$ such that $U\subseteq H$, and $xH$ supports finite
$H$-invariant Borel measure.
\end{thm}

Note that the proofs of Dani, Raghavan, Veech are different from the proofs that are presented
here. In fact, their proofs can be considered as the first important steps towards the general result on topological rigidity
-- Theorem \ref{th_rat}.

We start the proof of Theorem \ref{th_dr} with a simple lemma:

\begin{lem} \label{lem_ir}
Let $\{v_i:i=1,\ldots, s\}\subseteq\mathbb{R}^n$, $1\le s\le n-2$, be linearly independent 
vectors such that $\left<v_i:i=1,\ldots, s\right>$ contains no nonzero rational vectors.
Then there exists $v_{s+1}$ such that $v_i$, $i=1,\ldots, s+1$, are linearly independent, and
$\left<v_i:i=1,\ldots, s+1\right>$ contains no nonzero rational vectors.
\end{lem}

\begin{proof}
Let $V=\left<v_i:i=1,\ldots, s\right>$. Since $s\le n-2$, for any $v\in\mathbb{R}^n$,
$\left<V,v\right>$ is a proper subspace of $\mathbb{R}^n$. Therefore one can take a
vector $v_{s+1}$ outside $\cup_{v\in\mathbb{Q}^n}\left<V,v\right>$. If 
$v=\sum_{i=1}^{s+1} \alpha_i v_i$ is rational and nonzero for some $\alpha_i\in\mathbb{R}$,
then $\alpha_{s+1}\ne 0$, and $v_{s+1}\in \left<V,v\right>$. This is a contradiction.
Thus, $v_{s+1}$ is as required.
\end{proof}

\begin{proof}[Proof of Theorem \ref{th_dr}]
It is easy to see that
if the condition of the theorem is not satisfied, the orbit cannot be
dense. The hard part is to show that the above condition implies density.
By Lemma \ref{lem_ir}, we may assume that $l=n-1$.

Denote $G=\hbox{SL}(n,\mathbb{R})$, $\Gamma=\hbox{SL}(n,\mathbb{Z})$, and
$$
U_0=\left(
\begin{tabular}{c|c}
$E$ & $*$\\
\hline
$0$ & $1$
\end{tabular}
\right),
$$
where $E$ is the identity $(n-1)\times (n-1)$ matrix. Let $g_0\in G$ be such that $g_0v_i=e_i$ for
$i=1,\ldots, n-1$. Here $\{e_i:i=1,\ldots,n\}$ is the standard basis of $\mathbb{R}^n$.
Then the stabilizer of $v$ in $G$ is $U=g_0^{-1}U_0g_0$. Note that any nontrivial
$U$-invariant subspace of $\mathbb{R}^n$ is contained in $\left<v_i:i=1,\ldots, n-1\right>$.
Consider $U$-orbit $\Gamma U\subset \Gamma\backslash G$. By Ratner's theorem (Theorem \ref{th_rat}),
$\overline{\Gamma U}=\Gamma H$ where $H$ is a closed connected subgroup of $G$ containing $U$, and 
$H\cap \Gamma$ is a lattice in $H$. Moreover by \cite[Proposition 3.2]{sh91},
$H$ is the connected component of the smallest real algebraic $\mathbb{Q}$-subgroup containing $U$, and
the radical of $H$ is unipotent. Let $R$ be the radical of $H$. Since $R$ is defined
over $\mathbb{Q}$ and unipotent, the space $V^R$ of $R$-fixed vectors is nonzero
and defined over $\mathbb{Q}$. Also $V^R$ is $H$-invariant because $R$ is normal in $H$.
Thus, if $V^R\ne \mathbb{R}^n$, $V^R\subseteq\left<v_i:i=1,\ldots, n-1\right>$.
However, this contradicts our hypothesis on $v$. Therefore, $V^R=\mathbb{R}^n$
and $R=1$, i.e. $H$ is semisimple. 
We claim that $H=G$. To simplify notations, we work with the group $H_0\stackrel{def}{=}g_0Hg_0^{-1}$.
Let $\mathfrak h$ and $\mathfrak u$ be the Lie algebras of $H_0$ and $U_0$ respectively.
Since $U_0\subseteq H_0$,  
$$
\mathfrak u=\left<E_{in}:i=1,\ldots, n-1\right>\subseteq\mathfrak h.
$$
Here $E_{ij}$ denotes a matrix with $1$ at the place $(i,j)$ and $0$ elsewhere.
Using that the Killing form $k(x,y)=\hbox{Tr}(xy)$ for $x,y\in \mathfrak{sl}(n,\mathbb{R})$
is nondegenerate on $\mathfrak h$, one shows that the projection map from $\mathfrak h$ to the space
$\left<E_{ni}:i=1,\ldots,n-1\right>$ with respect to the basis $\{E_{ij}\}$ is surjective.
Thus for $i=1,\ldots,n-1$, there exists $h_i=E_{ni}+\tilde{h}_i\in \mathfrak h$ with 
$\tilde{h}_i$ in the normalizer of $\mathfrak{u}$. Then
$$
\mathfrak h\supseteq [h_i,\mathfrak u]+\mathfrak u\supseteq [E_{ni},\mathfrak u],\quad i=1,\ldots,n-1.
$$
It follows that $\mathfrak h=\mathfrak{sl}(n,\mathbb{R})$ and $H=G$. Thus, $\overline{\Gamma U}=G$.
Finally, 
\begin{equation}\label{eq_lll}
\overline{\Gamma v}= \overline{\Gamma Uv}\supseteq\overline{\Gamma U}v=Gv=\mathcal{V}_{n,l}.
\end{equation}
\end{proof}

\begin{proof}[Proof of Theorem \ref{th_ve}]
It is sufficient to prove the claim for $l=n-1$.
 
Let $U$ be as in the proof of Theorem \ref{th_dr}.
By (\ref{eq_lll}), we just need to show that $\Gamma U$ is dense in $G$.
By Ratner's theorem (Theorem \ref{th_rat}),
$\overline{\Gamma U}=\Gamma H$ where $H$ is 
a closed connected subgroup of $G$ containing $U$, and  $H\cap \Gamma$ is a lattice in $H$.
By Lemma 3.8 and Proposition 3.10 from \cite{sh91}, one of the following two possibilities holds:
$H$ is reductive, or $W\cap\Gamma$ is a lattice in $W$ where $W$ is the unipotent
radical of a proper parabolic subgroup of $G$. 
Since $\Gamma$ is cocompact,
it follows from Godement's criterion
that $\Gamma$ has no nontrivial unipotent elements. This contradicts the second
possibility. Thus, $H$ is reductive, and the Killing form is nondegenerate on
the Lie algebra of $H$.
Now one can show by the same argument as in the proof of Theorem \ref{th_dr} that $H=G$.
Hence, $\Gamma U$ is dense in $\Gamma\backslash G$.
This implies Theorem \ref{th_ve}.
\end{proof}

\subsection{Iwasawa decomposition for $\hbox{\rm SL}(n,\mathbb{R})$} \label{sec_iwas}

\hspace{1cm}\vspace{0.2cm}

Fix $l=1,\ldots,n-1$.

For $s=(s_1,\ldots,s_{n})\in\mathbb{R}^{n}$, $\sum_{i=1}^n s_i=0$,
define
$$
a(s)=\hbox{diag}(e^{s_1},\ldots, e^{s_n})\in\hbox{SL}(n,\mathbb{R}).
$$
For a vector $s$ as above, define decomposition 
$$
s=s^{_-}+s^{_+}
$$
with $s^{_-}=(s_1,\ldots,s_l,r,\ldots,r)$, $r=\frac{1}{n-l}(-s_1-\cdots-s_l)$,
$s^{_+} = s-s^{_-}$.
Note that $r$ is chosen so that $a(s^{_+}),a(s^{_-})\in\hbox{SL}(n,\mathbb{R})$.

For $t=(t_{ij}:1\le i<j\le l)$, $t_{ij}\in\mathbb{R}$, denote by $n^{_-}(t)$ the unipotent 
upper triangular matrix which entries above diagonal are equal $t_{ij}$ for $i<j\le l$ and $0$ otherwise.
Similarly, for $t=(t_{ij}:1\le i<j\le n, j>l)$, $t_{ij}\in\mathbb{R}$, denote by $n^{_+}(t)$ the unipotent 
upper triangular matrix which entries above diagonal are equal $t_{ij}$ for $1\le i<j\le n$, $j>l$ and $0$ otherwise.

We use the following notations:
\begin{eqnarray}
G&=&\hbox{SL}(n,\mathbb{R}),\nonumber\\
K&=&\hbox{SO}(n,\mathbb{R}),\nonumber\\
A^o_{l-}&=&\left\{ a(s^{_-})\;:\;s\in\mathbb{R}^n,\; \sum_{i=1}^n s_i=0\right\}\nonumber\\
A^o_{l+}&=&\left\{ a(s^{_+})\;:\; s\in\mathbb{R}^n,\; \sum_{i=1}^n s_i=0\right\},\nonumber\\
A^o&=&A^o_{l-}A^o_{l+}=\left\{ a(s)\;:\; s\in\mathbb{R}^n,\; \sum_{i=1}^n s_i=0\right\},\nonumber\\
N_{l-}&=&\{n^{_-}(t)\;:\; t_{ij}\in\mathbb{R},\; 1\le i<j\le l\},\nonumber\\
N_{l+}&=&\{n^{_+}(t)\;:\; t_{ij}\in\mathbb{R},\; 1\le i<j\le n, j>l\},\nonumber\\
N&=&N_{l-}N_{l+}=\hbox{``unipotent upper triangular group''},\nonumber\\
B^o_l&=&A^o_{l+}N_{l+}=N_{l+}A^o_{l+}. \nonumber
\end{eqnarray}

Denote by $dk$ the normalized Haar measure on $K$.

Let 
$$
dn^{_+}=dt^{+}=\prod_{i<j\le l}dt_{ij}\quad\hbox{and}\quad dn^{_-}=dt^{-}=\prod_{\max (i,l)<j}dt_{ij}.
$$
These measures are Haar measures for $N_{l+}$ and $N_{l-}$ respectively.
The subgroup $N_{l+}$ is normal in $N$, and the product map 
$$
N_{l-}\times N_{l+}\rightarrow N
$$
is a diffeomorphism.  Also the image of product of $dn^{_-}$ and $dn^{_+}$ under this map
is a Haar measure on $N$. Let us denote it by $dn$:
\begin{equation}\label{eq_dn}
\int_N f(n)dn=\int_{N_{l-}\times N_{l+}}f(n^{_-}n^{_+})dn^{_-}dn^{_+},\quad f\in C_c(N).
\end{equation}

Haar measures on $A^o_{l-}$ and $A^o_{l+}$ are defined by
$$
da^{_-}=ds^{_-}=\prod_{i=1}^l ds^{_-}_i\quad\hbox{and}\quad da^{_+}=ds^{_+}=\prod_{i=l+1}^{n-1} ds_i^{_+}
$$
respectively. Then a Haar measure $da$ on $A^o=A^o_{l-}A^o_{l+}$
is the product measure:
\begin{equation}\label{eq_da}
\int_{A^o}f(a)da=\int_{A^o_{l-}\times A^o_{l+}} f(a^{_-}a^{_+})da^{_-}da^{_+},\quad f\in C_c(A^o).
\end{equation}

The product map $A^o_{l+}\times N_{l+}\rightarrow B^o_l$ is a diffeomorphism, and
the image of the product measure under this map is a left Haar measure on $B^o_l$.
Denote this measure by $\lambda_l$. Then a right Haar measure $\varrho_l$ on $B^o_l$ can be defined by
\begin{equation}\label{eq_rho_l}
\varrho_l (f)=\int_{B^o_l} f(b^{-1}) \lambda_l(b)
=\int_{A^o_{l+}\times N_{l+}} f(a(s^{_+})n^{_+})\delta_l^+(s^{_+}) ds^{_+}dn^{_+} 
\end{equation}
for $f\in C_c(B^o_l)$, where
\begin{equation}\label{eq_dp}
\delta_l^+(a)=\delta_l^+(s)=\exp\left\{2\sum_{i=l+1}^n (n-i)s_i\right\}
\end{equation}
for $a=\hbox{diag}(e^{s_1},\ldots,e^{s_n})\in\hbox{SL}(n,\mathbb{R})$.

By Iwasawa decomposition, the map
\begin{equation} \label{eq_iwa0}
(k,n,a)\mapsto kna: K\times N\times A^o\rightarrow G
\end{equation}
is a diffeomorphism. 

\begin{lem} \label{l_vol}
Let $e^0=(e_1,\ldots, e_l)$ be the standard orthonormal frame in $\mathbb{R}^n$.
Then for $k\in K$, $n\in N$, and $a=\hbox{\rm diag}(e^{s_1},\ldots,e^{s_n})\in A^o$,
$$
\hbox{\rm \rm Vol}(knae^0)=\exp\left\{\sum_{i=1}^l s_i\right\}.
$$
\end{lem}

\begin{proof}
Let $g=kna$. Recall that Iwasawa decomposition is proved using Gramm-Schmidt
orthogonalization for basis $v_i=ge_i$. Let $\{w_i\}$ be an orthonormal basis
such that $\left<v_k:k\le i\right>=\left<w_k:k\le i\right>$ for $1\le i\le n$.
Then $e^{s_i}=v_i\cdot w_i$, i.e. $e^{s_i}$ is the length of projection of $v_i$
onto the orthogonal complement of $\left<v_k:k\le i-1\right>$ in
$\left<v_k:k\le i\right>$. Now the statement is obvious.
\end{proof}

Define
\begin{equation}\label{eq_dm}
\delta_l^-(a)=\delta_l^-(s)=\exp\left\{\sum_{i=1}^l s_i\right\},
\end{equation}
where $a=\hbox{diag}(e^{s_1},\ldots,e^{s_n})\in\hbox{SL}(n,\mathbb{R})$.

The image of the product measure under the map (\ref{eq_iwa0}) is
a Haar measure on $G$ \cite[Proposition X.1.12]{hel}. Let us denote this image by $\mu$:
\begin{equation}\label{eq_mu}
\int_G fd\mu=\int_{K\times N\times A^o}f(kna)dkdnda,\quad f\in C_c(G).
\end{equation}
For a lattice $\Gamma$ in $G$, there exists
a measure $\bar \mu$ on $\Gamma\backslash G$ such that
\begin{equation} \label{eq_mubar}
\int_G fd\mu=\int_{\Gamma\backslash G}\sum_{\gamma\in\Gamma} f(\gamma g)d\bar\mu(g),\quad f\in C_c(G).
\end{equation}

For our purposes, we modify the Iwasawa decomposition
as follows:
\begin{equation} \label{eq_iwasawa}
(k,n^{_-},a^{_-},b)\mapsto kn^{_-}a^{_-}b: K\times N_{l-}\times A^o_{l-}\times B^o_l\rightarrow G.
\end{equation}
Since $A^o_{l-}$ normalizes $B^o_{l}$, this map is a diffeomorphism too.

Fix $g_0\in G$. By (\ref{eq_iwasawa}) the map
\begin{equation} \label{eq_decomp}
(k,n^{_-},a^{_-},b)\mapsto kn^{_-}a^{_-}g_0b: K\times N_{l-}\times A^o_{l-}\times (B^o_l)^{g_0}\rightarrow G
\end{equation}
is a diffeomorphism. Here we use notation: $X^{g_0}=g_0^{-1}Xg_0$.

For $g\in G$ define $k_g\in K$, $a^{_\pm}_g\in A^o_{l\pm}$, and $n^{_\pm}_g\in N_{l\pm}$ such that
\begin{equation} \label{eq_decomp1}
g=k_gn^{_-}_ga^{_-}_gg_0 (a^{_+}_gn^{_+}_g)^{g_0}.
\end{equation}
Also define
$$
b^{-}_g=n^{_-}_ga^{_-}_g,\quad b_g=a^{_+}_gn^{_+}_g,
$$
and $b'_g\in\hbox{GL}(l,\mathbb{R})$ such that
\begin{equation} \label{eq_def_c}
b^{-}_g=\left(\begin{tabular}{r|r}
$b'_g$ & $0$\\
\hline
$0$ & $*$
\end{tabular}
\right).
\end{equation}

\begin{lem} \label{l_haar_measure}
For $f\in C_c(G)$ and $g_0\in G$,
$$
\int_G fd\mu=\mathop{\int}_{K\times N_{l-}\times A^o_{l-}\times B^o_{l}} f(kn^{_-}a(s^{_-})g_0b^{g_0})\delta_l^-(s^{_-})^n dkdn^{_-}ds^{_-}d\varrho_l(b),
$$
where $\delta_l^-$ is defined in (\ref{eq_dm}).
\end{lem}

\begin{proof}
By (\ref{eq_mu}), (\ref{eq_dn}), and (\ref{eq_da}),
\begin{equation}\label{eq_h1}
\int_G fd\mu=\mathop{\int}_{K\times N_{l-}\times N_{l+}\times A^o_{l-}\times A^o_{l+}} f(kn^{_-}n^{_+}a^{_-}a^{_+})dkdn^{_-}dn^{_+}da^{_-}da^{_+}.
\end{equation}
The Jacobian of the map 
$$
N_{l+}\rightarrow N_{l+}: n\mapsto a^{-1}n a
$$
for $a=a(s^{_-})a(s^{_+})$ is equal to $\delta_l^- (s^{_-})^{-n}\delta_l^+(s^{_+})^{-1}$. Thus,
it follows from (\ref{eq_rho_l}) and (\ref{eq_h1}) that
$$
\int_G fd\mu=\mathop{\int}_{K\times N_{l-}\times A^o_{l-}\times B^o_l} f(kn^{_-}a(s^{_-})b)\delta_l^- (s^{_-})^ndkdn^{_-}ds^{_-}d\varrho_l(b).
$$
Then since $G$ is unimodular,
\begin{eqnarray*}
\int_G fd\mu &=& \int_G f(g g_0) d \mu (g)\\
&=&\mathop{\int}_{K\times N_{l-}\times A^o_{l-}\times B^o_l} f(kn^{_-}a(s^{_-})bg_0)\delta_l^- (s^{_-})^ndkdn^{_-}ds^{_-}d\varrho_l(b)\\
&=&\mathop{\int}_{K\times N_{l-}\times A^o_{l-}\times B^o_l} f(kn^{_-}a(s^{_-})g_0b^{g_0})\delta_l^- (s^{_-})^ndkdn^{_-}ds^{_-}d\varrho_l(b).
\end{eqnarray*}
\end{proof}

\begin{lem} \label{l_frame_measure}
Let $e^0=(e_1,\ldots,e_l)$ be the standard frame in $\mathbb{R}^n$.
For $f\in C_c({\mathcal{V}}_{n,l})$,
\begin{equation} \label{eq_fm}
\int_{{\mathcal{V}}_{n,l}} f(v)dv=d_{n,l}\int_{K\times N_{l-}\times A^o_{l-}} f(kn^{_-}a(s^{_-})e^0) \delta_l^-(s^{_-})^n dkdn^{_-}ds^{_-},
\end{equation}
where $d_{n,l}$ is a constant computed in (\ref{eq_c0}).
\end{lem}

\begin{proof}
The measure on the left side of (\ref{eq_fm}) is $G$-invariant. We claim that the measure on the right
side is $G$-invariant too. It is easy to see that the map
$$
gB^o_l\mapsto ge^0:G/B^o_l\rightarrow \mathcal{V}_{n,l}
$$
is proper. Thus, every function $f\in C_c({\mathcal{V}}_{n,l})$ can be lifted to a function $f_1\in C_c(G/B^o_l)$.
Then the function $f$ can be represented as
$$
f(ge^0)=f_1(g B^o_l)=\int_{B_l^o}f_2(gb)d\varrho_l(b)
$$
for some $f_2\in C_c(G)$ (see \cite[Ch.~1]{rag}).
Then by Lemma \ref{l_haar_measure},
$$
\int_{K\times N_{l-}\times A^o_{l-}} f(kn^{_-}a(s^{_-})e^0) \delta_l^-(s^{_-})^n dkdn^{_-}ds^{_-} =\int_G f_2d\mu.
$$
It follows that the measure on the right side of (\ref{eq_fm}) is $G$-invariant. By uniqueness of Haar measure,
the integrals are equal up to a scalar multiple $d_{n,l}$. This constant
is computed in the Appendix.
\end{proof}

\subsection{Volume estimates}

\hspace{1cm}\vspace{0.2cm}

For a set $S\subseteq G$ and $T>0$, define
$$
S_T=\{s\in S: \|s\|<T\}.
$$
We compute the asymptotics of $\varrho_l (B^o_{l,T})$ as $T\rightarrow\infty$:

\begin{lem}\label{lem_BTasy}
\begin{equation}\label{eq_BTasy}
\varrho_l (B^o_{l,T})\sim \gamma_{n,l} T^{(n-1)(n-l)}\quad\hbox{as}\quad T\rightarrow\infty,
\end{equation}
where the constant $\gamma_{n,l}$ is given in (\ref{eq_gamnl}).
\end{lem}
The proof, which is given in the Appendix, follows the method of Duke, Rudnick, Sarnak \cite{drs}.

For $C\in\mathbb{R}$, define
\begin{eqnarray*}
A_{l+}^C &=& \{a(s^{_+}): s_i^{_+}>C, i=l+1,\ldots, n-1\},\\
B_l^C &=& A_{l+}^CN_{l+}.
\end{eqnarray*}
The following ``measure concentration'' result plays a crucial role in our proof.

\begin{lem}\label{lem_BTC2}
For $C\in\mathbb{R}$,
$$
\varrho_l(B^C_{l,T})\sim \varrho_l(B^o_{l,T})\;\;\hbox{as}\;\; T\rightarrow\infty.
$$
\end{lem}
This lemma is proved in the Appendix.

\subsection{Uniform distribution}

\hspace{1cm}\vspace{0.2cm}

\begin{thm}\label{th_erg_gen}
Let $\Gamma$ be a lattice in $G=\hbox{SL}(n,\mathbb{R})$. Fix $g_0\in (KN_{l-}A^o_{l-})^{-1}$ such that
$\Gamma (B^o_l)^{g_0}$ is dense in $G$. 
Let $Y=KN_{l-}A^o_{l-}g_0$, and $\nu_1$ be a measure on $Y$ such that
\begin{equation}\label{eq_nuBl}
\int_G fd\mu=\int_{Y}\int_{B^o_l} f(yb^{g_0})d\varrho_l(b)d\nu_1(y),\;\; f\in C_c(G).
\end{equation}
Let $\nu$ be a measure on $G/(B^o_l)^{g_0}$ induced by $\nu_1$.
For $T>0$ an $\Omega\subseteq G/(B^o_l)^{g_0}$, denote
\begin{equation}\label{eq_NTBl}
N_T(\Omega,g_0)=|\{\gamma\in\Gamma:\|\gamma\|<T,\gamma(B^o_l)^{g_0}\in\Omega\}|.
\end{equation}
Then for relatively compact Borel subset $\Omega$ of $G/(B^o_l)^{g_0}$ such that
$\nu(\partial\Omega)=0$,
$$
N_T(\Omega,g_0)\sim \varrho(B^o_{l,T})\frac{\delta_l^-(a_0)^{1-n}}{\bar\mu(\Gamma\backslash G)}\int_\Omega \frac{d\nu(x)}{\delta_l^-(a_x^{_-})}
\quad\hbox{as}\quad T\rightarrow\infty,
$$
where $a_0$ and $a_x^{_-}$ are the $A_{l-}^o$-components of $g_0^{-1}$ and $x$ with respect to decomposition
(\ref{eq_iwasawa}) respectively.
\end{thm}

Note that a similar result holds for every $g_0\in G$ because of the decomposition (\ref{eq_iwasawa}).

\begin{proof}
Write $g_0^{-1}=k_0b^{-}_0$ for $k_0\in K$, $b^{-}_0=n_0a_0\in N_{l-}A^o_{l-}$.

It is convenient to use decomposition (\ref{eq_decomp1}).
The product map $Y\times (B^o_l)^{g_0}\rightarrow G$ is a diffeomorphism.
For $g\in G$, denote $y_g=k_gn^{_-}_g a^{_-}_g g_0$, the $Y$-component of $g$. The map
$$
\alpha:Y\rightarrow G/(B^o_l)^{g_0}:y\mapsto y(B^o_l)^{g_0}
$$
is a diffeomorphism.
Clearly, $\gamma (B^o_{l})^{g_0}\in \Omega$ iff $y_\gamma\in \alpha^{-1} (\Omega)$.

For $g\in G$,
\begin{equation}\label{eq_nm1}
\|g\|=\|k_g b^{-}_g b_g g_0\|=\|b^-_g b_g (b^{-}_0)^{-1}\|.
\end{equation}
Note that for $\left(\begin{tabular}{c|c}
$*$ & $X$\\
\hline
$0$ & $Y$
\end{tabular}
\right)\in\hbox{SL}(n,\mathbb{R})$,
\begin{equation}\label{eq_nm2}
b^-_g\left(\begin{tabular}{c|c}
$*$ & $X$\\
\hline
$0$ & $Y$
\end{tabular}
\right)
(b^{-}_0)^{-1}=
\left(\begin{tabular}{c|c}
$*$ & $\beta_0^{-1} (b_g'\cdot X)$\\
\hline
$0$ & $\beta_0^{-1}\beta_gY$
\end{tabular}
\right),
\end{equation}
where
\begin{equation}\label{eq_beta}
\beta_0=\det (b_0')^{-\frac{1}{n-l}}\quad\hbox{and}\quad \beta_g=\det (b_g')^{-\frac{1}{n-l}}
\end{equation}
(here $b_0'$ and $b_g'$ are defined as in (\ref{eq_def_c})).
Put
\begin{equation}\label{eq_new2}
c_g=b^-_g\quad\hbox{and}\quad m_g=|\beta_0^{-1}\beta_g|.
\end{equation}
It follows from (\ref{eq_nm1}) and (\ref{eq_nm2}) that for $g\in G$,
$$
\|g\|^2=m_g^2\|\hat c_g(b_g)\|^2+e_g,
$$
where $e_g$ is a continuous function depending only on the $b^-$-components of $g$ with respect
to decomposition (\ref{eq_decomp1}).
 
Using previous notations, we have
\begin{equation} \label{eq_nt}
N_T(\Omega,g_0)=|\{\gamma\in\Gamma:m_\gamma^2\|\hat c_\gamma (b_\gamma)\|^2+e_\gamma<T^2,y_\gamma\in \alpha^{-1} (\Omega)\}|.
\end{equation}
To derive asymptotics of $N_T(\Omega,g_0)$, we can use Proposition \ref{pro_assym} with $H=B^o_l$, $h_\gamma=b_\gamma$, and $d_\gamma=e$.
Note that by Lemma \ref{lem_BTasy} the condition (\ref{eq_H1}) for $H=B^o_l$ is satisfied with
$M=(n-1)(n-l)$, and by Theorem \ref{th_ergodic} below, the condition (\ref{eq_H2}) holds.
Therefore, applying Proposition \ref{pro_assym}, we get
$$
N_T(\Omega,g_0)\sim \frac{\varrho(B^o_{l,T})}{\bar\mu(\Gamma\backslash G)}\int_{\alpha^{-1}(\Omega)}\frac{1}{m_y^{(n-1)(n-l)}\cdot\Delta_H(c_y)}d\nu_1(y)
$$
as $T\rightarrow\infty$, where $\Delta_H$ is defined in (\ref{eq_DH}). By (\ref{eq_beta}) and (\ref{eq_new2}),
$$
m_g=(\delta_l^-(a_0)^{-1}\delta_l^-(a_g^{_-}))^{-\frac{1}{n-l}}.
$$
Also
$$
\Delta_{B^o_l}(c_g)=
\frac{\det (b_g')^{n-l}}{\beta_g^{l(n-l)}}
=\delta_l^-(a_g^{_-})^{n}.
$$
Thus,
\begin{eqnarray*}
\int_{\alpha^{-1}(\Omega)}\frac{1}{m_y^{(n-1)(n-l)}\cdot\Delta_H(c_y)}d\nu_1(y)
&=&\int_{\alpha^{-1}(\Omega)}\frac{\delta_l^-(a_0)^{1-n}}{\delta_l^-(a_y^{_-})}d\nu_1(y)\\
&=&\delta_l^-(a_0)^{1-n}\int_\Omega \frac{d\nu(x)}{\delta_l^-(a^{_-}_x)}.
\end{eqnarray*}
This proves the theorem.
\end{proof}

\begin{proof}[Proof of Theorem \ref{th_frames00}]
For some $g_0\in (KN_{l-}A^o_{l-})^{-1}$, $v^0=g_0^{-1}e^0$ where $e^0=(e_1,\ldots,e_l)$ is the standard frame.
The condition that $\Gamma v^o$ is dense in $\mathcal{V}_{n,l}$ is equivalent to $\Gamma G_l^{g_0}$
is dense in $G$ where
$$
G_l=\left(\begin{tabular}{c|c}
$E$ & $*$\\
\hline
$0$ & $\hbox{SL}(n-l,\mathbb{R})$
\end{tabular}
\right).
$$
Since $B^o_l$ is epimorphic in $G_l$, it follows from \cite[Corollary 1.3]{sw} that $\Gamma (B^o_l)^{g_0}$ is dense in $G$.
  
Consider a map
$$
\alpha:G/(B^o_l)^{g_0}\rightarrow \mathcal{V}_{n,l}\simeq G/(G_l)^{g_0}: g(B^o_l)^{g_0}\mapsto gv^0.
$$
Note that this map is proper and $G$-equivariant.
Put $\Omega^*=\alpha^{-1}(\Omega)$. Then $\Omega^*$ is relatively compact, and
$N_T(\Omega,v^0)=N_T(\Omega^*,g_0)$, where $N_T(\Omega^*,g_0)$ is defined in (\ref{eq_NTBl}).

Let $\nu$ be the measure on $G/(B^o_l)^{g_0}$ defined in (\ref{eq_nuBl}).
It follows from Lemmas \ref{l_haar_measure} and \ref{l_frame_measure} that $\alpha(\nu)=d_{n,l}^{-1}dv$,
where $d_{n,l}$ is defined in (\ref{eq_c0}) and $dv$ is the Lebesgue measure on $\mathcal{V}_{n,l}$.

One can check that $\alpha(\partial\Omega^*)\subseteq\partial\Omega$.
Therefore,
$$
\nu(\partial\Omega^*)\le\nu(\alpha^{-1}(\partial\Omega))=d_{n,l}^{-1}\int_{\partial\Omega} dv=0.
$$
By Theorem \ref{th_erg_gen},
$$
N_T(\Omega,v^0)\sim \varrho(B^o_{l,T})\frac{\delta_l^-(a_0)^{1-n}}{\bar\mu(\Gamma\backslash G)}\int_{\Omega^*} \frac{d\nu(x)}{\delta_l^-(a_x^{_-})}\quad\hbox{as}\quad T\rightarrow\infty.
$$
By Lemma \ref{l_vol}, $\delta_l^-(a_0)=\hbox{\rm Vol}(v^0)$. Using decomposition (\ref{eq_decomp1}), we have
$gv^0=k_gn_g^{_-}a_g^{_-}e^0$. Thus, by Lemma \ref{l_vol}, $\delta_l^-(a^{_-}_g)=\hbox{\rm Vol}(gv^0)$.
Hence,
$$
N_T(\Omega,v^0)\sim \varrho(B^o_{l,T})\frac{\hbox{\rm Vol}(v^0)^{1-n}}{d_{n,l}\bar\mu(\Gamma\backslash G)}\left(\int_{\Omega} \frac{dv}{\hbox{\rm Vol}(v)}\right) \quad\hbox{as}\quad T\rightarrow\infty.
$$
Finally, using Lemma \ref{lem_BTasy} and (\ref{eq_c0}), we have
\begin{equation}\label{eq_NTasy}
N_T(\Omega,v^0)\sim a_{n,l}\frac{\hbox{\rm Vol}(v^0)^{1-n}}{\bar\mu(\Gamma\backslash G)}\left(\int_{\Omega} \frac{dv}{\hbox{\rm Vol}(v)}\right) T^{(n-1)(n-l)}\quad\hbox{as}\quad T\rightarrow\infty,
\end{equation}
where
\begin{equation}\label{eq_anl}
a_{n,l}=\frac{\gamma_{n,l}}{d_{n,l}}.
\end{equation}
The constants $\gamma_{n,l}$ and $d_{n,l}$ are computed in the Appendix.
\end{proof}

The proof of Corollary \ref{th_frames} is presented in the Appendix.

\section{Ergodic Theorem}\label{sec_erg}

The main result of this section is the following ergodic theorem along 
balls in $B^o_l$.

\begin{thm} \label{th_ergodic}
Let $1\le l\le n-1$.
Let $\Gamma$ be a lattice in $G$, and $y\in \Gamma\backslash G$ be such that $yB_l^o$ is dense in $\Gamma\backslash G$.
Denote by $\nu$ the probability $G$-invariant measure on $\Gamma\backslash G$.
Then for any $\tilde{f}\in C_c(\Gamma\backslash G)$,
$$
\frac{1}{\varrho_l(B^o_{l,T})}\int_{B^o_{l,T}} \tilde{f}(yb^{-1})d\varrho_l(b)\rightarrow \int_{\Gamma\backslash G} \tilde{f}d\nu\quad\hbox{as}\quad T\rightarrow\infty.
$$
\end{thm}

If $l=n-1$, the group $B^o_l$ is unipotent, and Theorem \ref{th_ergodic} is a special case of the result
of Shah \cite{sh94}. Thus, we may assume that $l<n-1$.

\subsection{Representations of $\hbox{\rm SL}(n,\mathbb{R})$}

\hspace{1cm}\vspace{0.2cm}

Before starting the proof, we prepare some auxiliary results on representations of
$G=\hbox{SL}(n,\mathbb{R})$.

Denote by $\mathfrak g$ the Lie algebra of $G$. We have the root space decomposition of $\mathfrak g_\mathbb{C}=\mathfrak g\otimes \mathbb{C}$:
$$
\mathfrak g_\mathbb{C}=\mathfrak g_0\oplus\sum_{i\ne j}\mathbb{C}E_{ij},
$$
where $\mathfrak g_0$ is the diagonal subalgebra of $\mathfrak g_\mathbb{C}$, and $E_{ij}$
is the matrix with 1 in position $(i,j)$ and $0$'s elsewhere.
It is convenient to identify $\mathfrak g_0$ with the space of vectors
$s=(s_1,\ldots,s_n)\in \mathbb{C}^n$, $\sum_i s_i=0$. Introduce the {\it roots} of $\mathfrak g_\mathbb{C}$:
$$
\alpha_{ij}(s)=s_i-s_j,\quad i\ne j,
$$
and the {\it fundamental weights} of $\mathfrak g_\mathbb{C}$:
$$
\omega_i(s)=s_1+\cdots+s_i,\quad 1\le i\le n-1.
$$
The {\it simple roots} of $\mathfrak g_\mathbb{C}$ are 
$\alpha_i=\alpha_{i,i+1}$, $i=1,\ldots,n-1$. For $i<j$,
$$
\alpha_{ij}=\alpha_i+\cdots+\alpha_j.
$$
The {\it dominant weights} are linear combinations with nonnegative integer coefficients of
the fundamental weights.
A {\it highest weight} of a finite-dimensional representation of $\mathfrak g_\mathbb{C}$
is a weight that is maximal with respect to the ordering on the dual space of $\mathfrak g_0$.
This weight is unique.
Irreducible representations of $\mathfrak g_\mathbb{C}$ are in one-to-one correspondence with
the dominant weights. The corresponding dominant weight is the highest weight of the representation.

Let $\mathfrak g_0^+$ be the Lie subalgebra of $\mathfrak g_\mathbb{C}$ that corresponds to
$A_{l+}^o$. That is, $\mathfrak g_0^+$ consists of diagonal matrices with entries
$$
(0,\ldots,0,s_{l+1},\ldots,s_n),\quad s_{l+1}+\cdots+s_n=0,
$$
on the diagonal.

\begin{lem}\label{l_dom}
Let $\pi$ be a representation of $G$ on a finite-dimensional complex vector space $V$.
Let 
\begin{equation}\label{eq_v0}
V_0=\{v\in V: \pi(B_l^o)v=v\}.
\end{equation}
Then every vector $\bar v\in V/V_0-\{0\}$ such that $\pi(N_{l+})\bar v=\bar v$ is a sum
of weight vectors of $\mathfrak g^+_0$ with nonzero dominant weights.
\end{lem}

\begin{proof}
First, we show that there are no nonzero vectors in $V/V_0$ fixed by $B^o_{l}$.
Let $W$ be the maximal $B_l^o$-invariant subspace on which $B_l^o$ acts unipotently.
The space $W$ can be constructed inductively as follows. Let $W_0$ be the space of
vectors fixed by $B^o_l$, $W_1\supseteq W_0$ be the space such that $W_1/W_0$
is the space of vectors in $V/W_0$ fixed by $B_l^o$, and so on. After finitely many
steps, we get $W$. We claim that $W=V_0$. Note $A^o_{l+}$ acts trivially on $W$.
Take $w\in W$. Suppose that $\tilde{\pi}(E_{ij})w\ne 0$ for some $E_{ij}\in\mathfrak n^{_+}$,
where $\mathfrak n^+$ is the Lie algebra of $N_{l+}$.
Then it is a weight vector with weight $\alpha_{ij}|_{\mathfrak g_0^+}\ne 0$
with respect to $\mathfrak g_0^+$. This is a contradiction. Thus,
$\tilde{\pi}(\mathfrak n^{_+})w=0$ for every $w\in W$, and $W=V_0$.
Clearly, the space $V/W$ does not contain any vectors that are fixed by $B^o_{l+}$.
This proves the claim.

Let 
$$
\tilde{G}=\left(
\begin{tabular}{c|c}
$E$ & $0$\\
\hline
$0$ & $\hbox{SL}(n-l,\mathbb{R})$
\end{tabular}
\right)\;\;\hbox{and}\;\;
\tilde{N}=\tilde{G}\cap N_{l+}.
$$
Since $A^o_{l+}\tilde{N}$ is an epimorphic subgroup of $\tilde{G}$, the space $V_0$ is $\tilde{G}$-invariant. 
Take a vector $\bar v\in V/V_0-\{0\}$ such that $\pi(N_{l+})\bar v=\bar v$.
By \cite[Lemma 15]{go} applied to $\tilde{G}$, $\bar v=\sum_{k=1}^m \bar v_k$ where $\bar v_k$,
$k=1,\ldots,m$, are weight vectors with
dominant weights $\lambda_k$ with respect to $\mathfrak g_0^+$.
Without loss of generality, $\lambda_k\ne\lambda_j$ for $k\ne j$.
For $E_{ij}\in\mathfrak n^+$, we have $\sum_{k=1}^m \tilde\pi(E_{ij})\bar v_k=0$.
Using that $\tilde\pi(E_{ij})\bar v_k$ is either $0$ or has weight $\lambda_k+\alpha_{ij}$,
we conclude that $\tilde\pi(E_{ij})\bar v_k=0$ for every $k=1,\ldots,m$.
Thus, $\pi(N_{l+})\bar v_k=\bar v_k$. Since $V_0=W$, the vector $\bar v_k$
cannot be fixed by $\pi(A_{l+}^o)$. Hence, $\lambda_k\ne 0$.
\end{proof}

We now modify slightly our notations. For $s=(s_{l+1},\ldots,s_n)\in\mathbb{R}^{n-l}$, $\sum_{i=l+1}^n s_i=0$, denote
$$
a^{_+}(s)=\hbox{diag}(1,\ldots,1,e^{s_{l+1}},\ldots,e^{s_n}).
$$

For $\beta>0$, define
\begin{equation}\label{eq_Db}
D(\beta)=\left\{t=(t_{ij}:\max(i,l)<j): \sum_{\max (i,l)<j}t_{ij}^2<\beta^2\right\}.
\end{equation}

\begin{lem}\label{lem_blow}
Let $\pi$ be a representation of $G$ on a finite-dimensional real vector space $V$, $V_0$ be defined as in (\ref{eq_v0}),
and $\bar V=V/V_0$. Introduce a norm on $\bar V$.
Then for every relatively compact subset $K\subseteq \bar V$ and $r>0$, there exists $\alpha\in (0,1)$
and $C_0>0$ such that for every $s$ such that $a^{_+}(s)\in A_{l+}^{C_0}$ and $x\in \bar V$ such that
$\|x\|>r$,
\begin{equation}\label{eq_sk}
\pi(a^{_+}(s)n^{_+}(D(e^{-\alpha s_{l+1}})))x\nsubseteq K.
\end{equation}
\end{lem}

\begin{proof}
The proof is the same as the proof of Lemma 16 in \cite{go}.
We will just sketch the main idea.

We need to get a lower estimate for 
$$
\sup\{\|\pi(a^{_+}(s))\pi(n^{_+}(t))x\|:t\in D(e^{-\alpha s_{l+1}})\}.
$$
Let $\bar W=\{\bar v\in \bar V:\pi(N_{l+})\bar v=\bar v\}$, and $\hbox{pr}_{\bar W}:\bar V\rightarrow\bar W$
is a projection on $\bar W$ that commutes with $\pi(a^{_+}(s))$.
By Lemma \ref{l_dom}, $\bar W$ is spanned by weight vectors of $A^o_{l+}$ with
nonzero dominant weights. 
Using that the character $s\mapsto e^{s_{l+1}}$ is the smallest nontrivial dominant weight
of $A^o_{l+}$, one concludes that for every $y\in\bar V$ and $a^{_+}(s)\in A_{l+}^C$ with $C>0$,\footnote{$A\ll B$ means $A < c\cdot B$ for some absolute constant $c>0$.}
\begin{equation}\label{eq_sk1}
\|\pi(a^{_+}(s))y\|\gg \|\pi(a^{_+}(s))\hbox{pr}_{\bar W}(y)\|\gg e^{s_{l+1}}\|\hbox{pr}_{\bar W}(y)\|.
\end{equation}
By a lemma due to Shah \cite{sh96} (see \cite[Lemma 13]{go}),
\begin{equation}\label{eq_sk2}
\sup\{\|\hbox{pr}_{\bar W}(\pi(n^{_+}(t))x)\|:t\in D(e^{-\alpha s_{l+1}})\}\gg (e^{-\alpha s_{l+1}})^d\|x\|.
\end{equation}
for some positive integer $d$. Combining (\ref{eq_sk1}) and (\ref{eq_sk2}), we get
\begin{equation}\label{eq_sk3}
\sup\{\|\pi(a^{_+}(s))\pi(n^{_+}(t))x\|:t\in D(e^{-\alpha s_{l+1}})\}\gg e^{cs_{l+1}}\|x\|,
\end{equation}
where $c=1-\alpha d$. Choose $\alpha\in (0,1)$ such that $c>0$. Then for $a^{_+}(s)\in A_{l+}^{C}$
the right hand side of (\ref{eq_sk3}) gets arbitrarily large as $C\rightarrow\infty$.
This proves (\ref{eq_sk}).
\end{proof}

\subsection{Proof of Theorem \ref{th_ergodic}}

\hspace{1cm}\vspace{0.2cm}

Now we are ready to start the proof of Theorem \ref{th_ergodic}.
Let $\mathcal{X}=(\Gamma\backslash G)\cup\{\infty\}$ be the one-point compactification of
$\Gamma\backslash G$. For $T>0$, define a normalized measure on $\mathcal{X}$ by
$$
\nu_T(\tilde{f})=\frac{1}{\varrho_l(B^o_{l,T})}\int_{B^o_{l,T}} \tilde{f}(yb^{-1})d\varrho_l(b),\quad \tilde{f}\in C_c(\Gamma\backslash G).
$$
We need to show that $\nu_T\rightarrow\nu$ as $T\rightarrow\infty$ in weak${^*}$ topology.
Since the space of normalized measures on $\mathcal{Z}$ is compact in weak${^*}$ topology,
it is enough to show that every limit point of $\nu_T$, $T\rightarrow\infty$,
is equal to $\nu$. Let $\nu_{T_i}\rightarrow\eta$ as $T_i\rightarrow\infty$
for some normalized measure $\eta$ on $\mathcal{X}$.
By Lemma \ref{lem_BTC2}, for every $C\in\mathbb{R}$,
\begin{equation}
\eta(\tilde{f})=\lim_{T_i\rightarrow\infty}\frac{1}{\varrho_l(B^o_{l,T_i})}\int_{B^C_{l,T_i}} \tilde{f}(yb^{-1})d\varrho_l(b)
\end{equation}

Let $U=\{n^{_+}(t)\in N:\; t_{ij}=0\;\hbox{for}\; i<j<n\}$.

\begin{lem}\label{l_u_inv}
The measure $\eta$ is $U$-invariant.
\end{lem}

Up to minor modifications, the proof is the same as the proof of Lemma 18 in \cite{go}.

\begin{lem}\label{l_a_tilde}
For $\alpha\in (0,1)$, define
\begin{equation}\label{eq_tildeA}
\tilde{A}_{l+,T}^C=\left\{a^{_+}(s)\in A_{l+,T}^C:\; (T^2-N(s)-l)^{1/2}>\exp\left(\mathop{\max}_{l<i<n} \{s_i\} -\alpha s_{l+1} \right)\right\},
\end{equation}
where $N(s)$ is defined in (\ref{eq_Ns}), and
$$
\tilde{B}_{l,T}^C=\left(\tilde{A}_{l+,T}^C N_{l+}\right)\bigcap B^o_{l,T}.
$$
Then for every $C>0$,
$$
\eta(\tilde{f})=\lim_{T_i\rightarrow\infty}\frac{1}{\varrho_l(B^o_{l,T_i})}\int_{\tilde{B}^C_{l,T_i}} \tilde{f}(yb^{-1})d\varrho_l(b).
$$
\end{lem}

The proof is routine computation based on Lemma \ref{lem_BTC} in the Appendix.
See Lemma 19 in \cite{go} or the proof of Lemma \ref{lem_BTC2} above for a similar argument.

Write $y=\Gamma g_0$ for some $g_0\in G$.
Let $q_s(t)=n^{_+}(t)^{-1}a^{_+}(s)^{-1}$.

Next, we review some deep results on distribution of polynomial trajectories
due to Dani, Margulis, Shah, and Ratner. See \cite{kss01} and \cite[\S 19]{st} for more comprehensive exposition.
These results will be applied to the polynomial map
$$
t\mapsto g_0 q_s(t):\mathbb{R}^m\rightarrow\hbox{SL}(n,\mathbb{R}),
$$
where $m=\frac{1}{2}(n-l)(n+l-1)$.

Let $\mathfrak g$ be the Lie algebra of $G$, and $V_G=\oplus_{i=1}^{\dim g}\wedge^i\mathfrak g$.
Fix a norm in $V_G$. For every Lie subgroup $H$ of $G$ with Lie algebra $\mathfrak h$, take a unit
vector $p_H\in\wedge^{\dim\mathfrak h}\mathfrak h\subseteq V_G$. Also define
$$
X(H,U)=\{g\in G:gU\subseteq Hg\}.
$$
Denote by $\mathcal{H}_\Gamma$ the family of all proper closed connected subgroups $H$ of $G$
such that $\Gamma\cap H$ is a lattice in $H$, and $\hbox{Ad}(\Gamma\cap H)$ is Zariski dense
in $\hbox{Ad}(H)$.

The {\it singular set} of $U$ is 
$$
Y=\bigcup_{H\in\mathcal{H}_\Gamma} \Gamma X(H,U)\subseteq\Gamma\backslash G.
$$
The set $Y$ is precisely the set of $y\in \Gamma\backslash G$ such that $yU$ is not dense
in $\Gamma\backslash G$.

The following facts and results will be used in the sequel:

\begin{enumerate}
\item[(I)] The set $\mathcal{H}_\Gamma$ is countable.
\item[(II)] For every $H\in\mathcal{H}_\Gamma$, $\Gamma\cdot p_H$ is discrete in $V_G$.
Thus, $\Gamma N^1_G(H)$ is closed in $\Gamma\backslash G$, where
$N_G^1(H)=\{g\in G:g\cdot p_H=p_H\}$.
\item[(III)] 
Assume that $\Gamma$ is not cocompact.
Then there exist closed subgroups $U_i$, $i=1,\ldots, r$, such that each $U_i$ is
the unipotent radical of a parabolic subgroup, $\Gamma U_i$ is closed in $\Gamma\backslash G$,
and for every $\varepsilon,\delta>0$, there exists a compact set $C\subseteq\Gamma\backslash G$
such that for every bounded open convex subset $D\subseteq\mathbb{R}^m$, one of the following holds:
\begin{enumerate}
\item[1.] There exist $\gamma\in\Gamma$ and $i=1,\ldots,r$ such that
$$
\mathop{\sup}_{t\in D}\|q_s(t)^{-1}g_0^{-1}\gamma\cdot p_{U_i}\|\le \delta.
$$
\item[2.] $\omega(\{t\in D:\Gamma g_0q_s(t)\notin C\})<\varepsilon\omega(D)$, where $\omega$ is the Lebesgue measure
on $\mathbb {R}^m$.
\end{enumerate}
\item[(IV)] Let $\varepsilon>0$ and $H\in\mathcal{H}_\Gamma$.
For every compact set $C\subseteq \Gamma X(H,U)$, there exists a compact set $F\subseteq V_G$
such that for every neighborhood $\Phi$ of $F$ in $V_G$ there exists a neighborhood
$\Psi$ of $C$ in $\Gamma\backslash G$ such that for every bounded open convex set
$D\subseteq\mathbb{R}^m$, one of the following holds:
\begin{enumerate}
\item[1.] There exists $\gamma\in \Gamma$ such that $q_s(D)^{-1}g_0^{-1}\gamma\cdot p_H\subseteq \Phi$.
\item[2.] $\omega(\{t\in D:\Gamma g_0q_s(t)\in \Psi\})<\varepsilon\omega(D)$.
\end{enumerate}
\end{enumerate}

(I) is proved in \cite[Theorem 1.1]{r91a} and \cite[Theorem 2.1]{dm93}.
For the proof of (II), see \cite[Theorem 3.4]{dm93}. (III) is a special case of \cite[Theorems 2.1--2.2]{sh96}.
(IV) is based on \cite[Proposition 5.4]{sh94}. It is formulated in \cite{sh96}.

To simplify our notations, we put $V=V_G$. Let $V_0$ be defined as in (\ref{eq_v0}).

\begin{lem} \label{lem_v0}
For $H\in\mathcal{H}_\Gamma$, $g_0^{-1}\Gamma\cdot p_H\subseteq V-V_0$.
\end{lem}

\begin{proof}
Suppose that $g_0^{-1}\gamma\cdot p_H\in V_0$ for some $\gamma\in\Gamma$.
Then $(\gamma^{-1}g_0B_l^og_0^{-1}\gamma)\cdot p_H=p_H$. Thus,
$\gamma^{-1}g_0B_l^og_0^{-1}\gamma\subseteq N_G^1(H)$. By (II), $\Gamma N^1_G(H)$
is not dense in $\Gamma\backslash G$. It follows that $\Gamma g_0B_l^o$
is not dense in $\Gamma\backslash G$ too. This is a contradiction.
\end{proof}

In the case when $\Gamma$ is not cocompact, we prove the following lemma:

\begin{lem}\label{lem_infty}
$\eta(\{\infty\})=0$.
\end{lem}

\begin{proof}
Write
\begin{equation}\label{eq_v01}
V=V_0\oplus V_1,
\end{equation}
where $V_1$ is $A_{l+}^o$-invariant complement. 
For a vector $v\in V$, denote by $v_0\in V_0$ and $v_1\in V_1$ its components with respect
to the decomposition (\ref{eq_v01}).
Fix norms on $V_0$ and $V_1$. Define a norm
on $V$ by
$$
\|v\|=\max\{\|v_0\|,\|v_1\|\},\quad v_0\in V_0,\; v_1\in V_1.
$$
The space $V_1$ is naturally isomorphic with $V/V_0$. The norm on $V_1$ induces a norm
on $V/V_0$ through this isomorphism.

We use (III). Let $\varepsilon,\delta>0$. Let
$$
P=\bigcup_{i=1}^r g_0^{-1}\Gamma\cdot p_{U_i},
$$
and $P_1=\{p\in P:\|p_0\|>\delta\}$. For $p\in P_1$,
\begin{equation}\label{eq_inf0}
\|q_s(0)^{-1}p\|\ge\|p_0\|>\delta.
\end{equation}
By Lemma \ref{lem_v0}, $P\subseteq V-V_0$, and by (II), $P$ is discrete. Therefore,
there exists $r>0$ such that $\|p_1\|>r$ for $p\in P-P_1$. Since the factor-map 
$V\rightarrow V/V_0$ is continuous and $B^o_l$-equivariant, for some $M>0$,
\begin{equation}\label{eq_inf1}
\|q_s(t)^{-1}\cdot v\|\ge M \|q_s(t)^{-1}\cdot \bar v\|,\quad v\in V.
\end{equation}
Now we apply Lemma \ref{lem_blow} with $K=\{\bar v\in V/V_0:\|\bar v\|\le \frac{\delta}{M}\}$.
There exist $\alpha\in (0,1)$ and $C_0>0$ such that for every $s$ such that
$a^{_+}(s)\in A_{l+}^{C_0}$ and every $\bar v\in V/V_0$ such that $\|\bar v\|>r$,
$$
q_s(D(e^{-\alpha s_{l+1}}))^{-1}\cdot \bar v \nsubseteq K.
$$
In particular, this holds for $\bar v=\bar p$ with $p\in P-P_1$.
Thus, by (\ref{eq_inf1}),
\begin{equation}\label{eq_inf2}
\mathop{\sup}_{t\in D(e^{-\alpha s_{l+1}})} \|q_s(t)^{-1}\cdot p\|>\delta
\end{equation}
for $p\in P-P_1$. In fact, (\ref{eq_inf2}) holds for $p\in P_1$ because of (\ref{eq_inf0}).
Thus, the case (a) of (III) does not occur when $a^{_+}(s)\in A^{C_0}_{l+,T_i}$ and
$D$ is a bounded open convex subset such that $D\supseteq D(e^{-\alpha s_{l+1}})$.
It follows that for some compact set $C\subseteq \Gamma\backslash G$,
\begin{equation}\label{eq_mC}
\omega(\{t\in D:\Gamma g_0q_s(t)\notin C\})<\varepsilon\omega(D).
\end{equation}
when $a^{_+}(s)\in A^{C_0}_{l+,T_i}$ and $D\supseteq D(e^{-\alpha s_{l+1}})$.

We have
\begin{equation} \label{eq_tA}
\eta(\tilde{f})=\lim_{T_i\rightarrow\infty}
\frac{1}{\varrho_l(B^o_{l,T_i})}\int_{\tilde{A}^C_{l+,T_i}}\int_{D_{s,T_i}}\tilde{f}(\Gamma
g_0q_s(t))\delta_l^+(s)dt^+ds^{_+},\quad \tilde{f}\in C_c(\Gamma\backslash G),
\end{equation}
where
\begin{equation}\label{eq_DsT}
D_{s,T_i}=\left\{n(t)\in N: \sum_{i\le l; l<j} t_{ij}^2+ \sum_{l<i<j} e^{2s_i}t_{ij}^2<T_i^2-N(s)-l\right\},
\end{equation}
and $N(s)$ is defined in (\ref{eq_Ns}).
Note that $D_{s,T_i}$ contains $D(\beta)$, which is defined in (\ref{eq_Db}), for 
$$
\beta<(T_i^2-N(s)-l)^{1/2}\exp\left(-\mathop{\max}_{l+1\le i\le n-1} \{s_i\}\right).
$$
When $a^{_+}(s)\in\tilde{A}^{C_0}_{l+,T_i}$, the right hand side is greater then $e^{-\alpha s_{l+1}}$
(see (\ref{eq_tildeA})). Therefore, $D_{s,T_i}\supseteq D(e^{-\alpha s_{l+1}})$ when
$a^{_+}(s)\in\tilde{A}^{C_0}_{l+,T_i}$. By (\ref{eq_mC}),
\begin{equation} \label{eq_mC1}
\omega \left(\left\{t\in D_{s,T_i}: \Gamma g_0 q_s(t)\notin C\right\}\right)< \varepsilon
\omega (D_{s,T_i})\quad \hbox{for}\;\; a^{_+}(s)\in\tilde{A}^{C_0}_{l+,T_i}.
\end{equation}

Let $\chi_C$ be the characteristic function of the set $C$.
Take $\tilde{f}\in C_c(\Gamma\backslash G)$
such that $\chi_C\le \tilde{f}\le 1$. Then using (\ref{eq_tA}) and (\ref{eq_mC1}), we get
\begin{eqnarray*}
\eta(\hbox{supp}(\tilde{f}))&\ge& \lim_{T_i\rightarrow\infty}
\frac{1}{\varrho_l(B^o_{l,T_i})}\int_{\tilde{A}^{C_0}_{l+,T_i}}\int_{D_{s,T_i}}\chi_C(\Gamma g_0
q_s(t))\delta_l^+(s)dt^+ds^{_+}\\
&\ge&\lim_{T_i\rightarrow\infty}
\frac{1}{\varrho_l(B^o_{l,T_i})}\int_{\tilde{A}^{C_0}_{l+,T_i}}(1-\varepsilon)\omega(D_{s,T_i})\delta_l^+(s)ds^{_+}\\
&=&(1-\varepsilon)\lim_{T_i\rightarrow\infty} \frac{\varrho_l(\tilde{B}^{C_0}_{l,T_i})}{\varrho_l(B^o_{l,T_i})}=1-\varepsilon.
\end{eqnarray*}
Hence, $\eta(\{\infty\})\le\eta(\hbox{supp}(\tilde{f})^c)\le\varepsilon$ for every $\varepsilon>0$.
\end{proof}

\begin{lem}\label{lem_Y}
$\eta(Y)=0$.
\end{lem}

\begin{proof}
Since $\mathcal{H}_\Gamma$ is countable, it is enough to show that $\eta(\Gamma X(H,U))=0$
for every $H\in \mathcal{H}_\Gamma$. Moreover, it is enough to show that
$\eta(C)=0$ for every compact set $C\subseteq \Gamma X(H,U)$.

We use the notations from the proof of Lemma \ref{lem_infty}, in particular, decomposition (\ref{eq_v01}).

We apply (IV). Take $\varepsilon>0$. Let $F$ be a compact subset of $V$ as in (IV).
Take a relatively compact neighborhood $\Phi$ of $F$. Let $\Psi\supset C$ be as in (IV).
Denote $P=g_0^{-1}\Gamma\cdot p_H$ and $P_1=\{p\in P:\|p_0\|>\delta\}$ with
$\delta=\sup\{\|v_0\|: v\in \Phi\}$. Then (\ref{eq_inf0}) holds, so that
\begin{equation}\label{eq_p1}
q_s(0)\cdot p\notin \Phi.
\end{equation}
As in the proof of the previous lemma, there exists $r>0$ such that
$\|p_1\|>r$ for $p\in P-P_1$, and applying Lemma \ref{lem_blow},
we deduce that there exist $\alpha\in (0,1)$ and $C_0>0$ such that for every $s$ such that
$a^{_+}(s)\in A_{l+}^{C_0}$ and every $p\in P-P_1$,
\begin{equation}\label{eq_p2}
q_s(D(e^{-\alpha s_{l+1}}))^{-1}\cdot p \nsubseteq \Phi.
\end{equation}
By (\ref{eq_p1}), (\ref{eq_p2}) holds for every $p\in P$. Thus, case (a) of (IV) fails.
Therefore, case (b) holds:
\begin{equation}\label{eq_mCnew}
\omega(\{t\in D:\Gamma g_0q_s(t)\in \Psi\})<\varepsilon\omega(D),
\end{equation}
when $a^{_+}(s)\in A^{C_0}_{l+,T_i}$ and $D$ is an open convex set such that $D\supseteq D(e^{-\alpha s_{l+1}})$.
Recall that $D_{s,T_i}$ was defined in (\ref{eq_DsT}). It is easy to see from
(\ref{eq_DsT}) that $D_{s,T_i}\supseteq D(e^{-\alpha s_{l+1}})$ when
$a^{_+}(s)\in\tilde{A}^{C_0}_{l+,T_i}$. Thus, (\ref{eq_mCnew}) holds for $D=D_{s,T_i}$
with $a^{_+}(s)\in\tilde{A}^{C_0}_{l+,T_i}$.

Take a function $\tilde{f}\in C_c(\Gamma\backslash G)$ such that $\tilde{f}=1$ on $C$,
$\hbox{supp}(\tilde{f})\subseteq \Psi$, and $0\le \tilde{f}\le 1$.
Let $\chi_\Psi$ be the characteristic function of $\Psi$.
Then using (\ref{eq_tA}) and (\ref{eq_mCnew}) with $D=D_{s,T_i}$, we get
\begin{eqnarray*}
\eta(C)\le \lim_{T_i\rightarrow\infty}
\frac{1}{\varrho_l(B^o_{l,T_i})}\int_{\tilde{A}^{C_0}_{l+,T_i}}\int_{D_{s,T_i}}\chi_\Psi(\Gamma g_0
q_s(t))\delta_l^+(s)dt^+ds^{_+}\\
\le \lim_{T_i\rightarrow\infty}
\frac{1}{\varrho_l(B^o_{T_i})}\int_{\tilde{A}^{C_0}_{l+,T_i}}\varepsilon\omega(D_{s,T_i})\delta_l^+(s)ds^{_+}
=\varepsilon\lim_{T_i\rightarrow\infty} \frac{\varrho_l(\tilde{B}^{C_0}_{l,T_i})}{\varrho_l(B^o_{T_i})}=\varepsilon.
\end{eqnarray*}
This shows that $\eta(C)=0$. Hence, $\eta (Y)=0$.
\end{proof}

By Lemma \ref{l_u_inv}, the measure $\eta$ is $U$-invariant.
Consider the ergodic decomposition of $\eta$ into $U$-ergodic measures.
By Ratner's measure classification \cite{r91a}, an ergodic component of $\eta$ is
either $G$-invariant or supported on $Y\cup\{\infty\}$. By Lemmas
\ref{lem_infty} and \ref{lem_Y}, the set of ergodic components of the second type
has measure $0$. Therefore, $\eta$ is $G$-invariant, and $\eta=\nu$.
This proves Theorem \ref{th_ergodic}.

\section{Uniform distribution for a lattice in $\hbox{Sp}(n,\mathbb{R})$}\label{sec_sp}

\subsection{Density of orbits}

\hspace{1cm}\vspace{0.2cm}

\begin{proof}[Proof of Theorem \ref{th_dr2}]
Clearly, the condition is neccesary for denseness. Suppose that the condition
is satisfied.
Let $\{e_i:i=1,\ldots,2n\}$ be the standard basis of $\mathbb{R}^{2n}$, and $e=(e_1,\ldots, e_n)$.
Then $e$ is an isotropic frame, and by Witt's Theorem, the space of isotropic
$n$-frames is $Ge$. The stabilizer of $e$ in $G$ is
$$
U_0=\left\{\left(
\begin{tabular}{c|c}
$E$ & $X$\\
\hline
$0$ & $E$
\end{tabular}
\right): {}^t X=X\right\}.
$$
Let $g_0\in G$ be such that $g_0v=e$. Then the stabilizer of $v$ in $G$ is
$U=g_0^{-1}U_0g_0$. It is not hard to check that any $U_0$-invariant subspace
is either contained in $\left<e_1,\ldots, e_n\right>$ or contains it.
It follows that any $U$-invariant subspace is either contained
in $V_0\stackrel{def}{=}\left<v_1,\ldots, v_n\right>$ or contains it.

As in the proof of Theorem \ref{th_dr}, $\overline{\Gamma U}=\Gamma H$ where $H$ is 
the connected component of the smallest
real $\mathbb{Q}$-algebraic subgroup containing $U$, and the radical of $H$ is unipotent.
Let $R$ be the radical of $H$. The space of $R$-fixed vectors $V^R$ is defined over
$\mathbb{Q}$ and $H$-invariant. Since $R$ is unipotent, $V^R\ne 0$.
Thus, by the condition on $v$, $V^R\nsubseteq V_0$. Then $V_0\subseteq V^R$.
Suppose that $V^R\ne \mathbb{R}^{2n}$. Since $V_0^\perp=V_0$, $J|_{V^R}$ 
is degenerate. Then $0\ne\hbox{Rad} (J|_{V^R})\subseteq V_0$. This is a contradiction
because the space $\hbox{Rad} (J|_{V^R})$ is defined over $\mathbb{Q}$.
Hence, $V^R=\mathbb{R}^{2n}$, $R=1$, and $H$ is semisimple.

We claim that $H=G$. Let $H_0=g_0Hg_0^{-1}$. Denote by
$\mathfrak g$, $\mathfrak h$, and $\mathfrak u$ the Lie algebras of $G$, $H_0$, and $U_0$ respectively.
The Killing form on $\mathfrak g$ is defined by $k(x,y)=\hbox{Tr}(xy)$ for
$x,y\in\mathfrak g$. Since $H_0$ is semisimple, $k$ is nondegenerate on $\mathfrak h$.

Recall the root decomposition for $\mathfrak g$. A Cartan subalgebra of $\mathfrak g$ is
$$
\mathfrak a=\{\hbox{diag}(h_1,\ldots, h_n, -h_1,\ldots, -h_n):h_i\in\mathbb{R}\}.
$$
Let $\alpha_i(h)=h_i$ for $h\in\mathfrak a$. The root system of $\mathfrak g$ is
$$
\Delta=\{\alpha_i-\alpha_j,\pm(\alpha_p+\alpha_q):1\le i\ne j\le n,1\le p\le q\le n\},
$$
and the following root space decomposition holds:
\begin{eqnarray}
\mathfrak g_{\alpha_i-\alpha_j}&=&\left<E_{i,j}-E_{j+n, i+n}\right>,\nonumber\\
\mathfrak g_{\alpha_p+\alpha_q}&=&\left<E_{p, q+n}+E_{q, p+n}\right>,\nonumber\\
\mathfrak g_{-\alpha_p-\alpha_q}&=&\left<E_{p+n, q}+E_{q+n, p}\right>,\nonumber\\
\mathfrak g &=& \mathfrak a\oplus\sum_{\alpha\in\Delta} \mathfrak g_\alpha.\label{eq_dec}
\end{eqnarray}
Since $U_0\subseteq H_0$,
$$
\mathfrak u=\sum_{1\le p\le q\le n}\mathfrak g_{\alpha_p+\alpha_q}\subseteq\mathfrak h.
$$
Note that $k(\mathfrak g_\alpha,\mathfrak g_\beta)=0$ if $\alpha+\beta\ne 0$.
Since the Killing form $k$ is nondegenerate on $\mathfrak h$,
the projection map from $\mathfrak h$ to the space $\sum_{1\le p\le q\le n}\mathfrak g_{-\alpha_p-\alpha_q}$
with respect to the decomposition (\ref{eq_dec}) is surjective.
Thus for $1\le p\le q\le n$, there exists $h_{pq}=x_{pq}+\tilde{h}_{pq}\in\mathfrak h$
where $x_{pq}$ is a generator of the space $\mathfrak g_{-\alpha_p-\alpha_q}$, and $\tilde{h}_{pq}$
is in the normalizer of $\mathfrak u$. Then
$$
\mathfrak h\supseteq [h_{pq},\mathfrak u]+\mathfrak u\supseteq [x_{pq},u],\quad 1\le p\le q\le n.
$$
Using that $[\mathfrak g_\alpha,\mathfrak g_\beta]=\mathfrak g_{\alpha+\beta}$
for $\alpha,\beta,\alpha+\beta\in\Delta$, we conclude that 
$$
\sum_{\alpha\in\Delta}\mathfrak g_\alpha\subseteq\mathfrak h.
$$
It follows that $\mathfrak h=\mathfrak g$, and $H=G$.
Finally, 
$$
\overline{\Gamma v}= \overline{\Gamma Uv}\supseteq\overline{\Gamma U}v=Gv=\mathcal{W}_n.
$$
\end{proof}

\subsection{Iwasawa decomposition for $\hbox{Sp}(n,\mathbb{R})$} \label{sec_iwas1}

\hspace{1cm}\vspace{0.2cm}

Let $G=\hbox{Sp}(n,\mathbb{R})$. We use the following notations:
\begin{eqnarray*}
K&=&G\cap\hbox{SO}(2n,\mathbb{R}),\\
N_+&=&\left\{
\left(\begin{tabular}{ll}
$E$ & $N$\\
$0$ & $E$
\end{tabular}
\right):
{}^tN=N\right\},\\
N_-&=&\left\{
\left(\begin{tabular}{cc}
$M$ & $0$\\
$0$ & ${}^tM^{-1}$
\end{tabular}
\right): M\;\hbox{is upper triangular, unipotent}\right\},\\
A&=&\left\{
\left(\begin{tabular}{cc}
$A$ & $0$\\
$0$ & $A^{-1}$
\end{tabular}
\right): A\;\hbox{is positive, diagonal}\right\},\\
B&=&N_-A,\\
N&=&N_-N_+.
\end{eqnarray*}
We have Iwasawa decomposition:
$$
(k,n,a)\mapsto kna:K\times N\times A\rightarrow G
$$
(see, for example, \cite[p.~286]{ter2}). This map is a diffeomorphism.
It is easy to check that the product map $N_-\times N_+\rightarrow N$ is a diffeomorphism,
$N_-$ normalizes $N_+$, and $A$ normalizes $N_+$. Thus, modified Iwasawa decomposition holds:
\begin{equation}\label{eq_spiw}
(k,b,n)\mapsto kbn:K\times B\times N_+\rightarrow G.
\end{equation}

Fix $g_0\in G$. We also have decomposition:
\begin{equation}\label{eq_spiw2}
(k,b,n)\mapsto kbg_0 n^{g_0}:K\times B\times N_+\rightarrow G.
\end{equation}
For $g\in G$, define $k_g\in K$, $b_g\in B$, and $n_g\in N_+$ such that
$$
g=k_gb_gg_0n_g^{g_0}.
$$
Also define $b_g'\in\hbox{GL}(n,\mathbb{R})$ such that
$$
b_g=\left(\begin{tabular}{cc}
$b_g'$ & $0$\\
$0$ & $^t (b_g')^{-1}$
\end{tabular}
\right).
$$

Let $\mu$ be a Haar measure on $G$, and $\bar\mu$ be the measure on $\Gamma\backslash G$ such that
$$
\int_G fd\mu=\int_{\Gamma\backslash G}\sum_{\gamma\in \Gamma} f(\gamma g)d\bar\mu(g),\quad f\in C_c(G).
$$
Let $\varrho$ be the Lebesgue measure on $N_+\simeq \mathbb{R}^\frac{n(n+1)}{2}$, and $\nu$
be the measure on $G/N_+^{g_0}$ such that
$$
\int_G fd\mu=\int_{G/N_+^{g_0}}\int_{N_+} f(gn^{g_0})d\varrho(n)d\nu(g),\quad f\in C_c(G).
$$
Note that
\begin{equation}\label{eq_last}
\varrho(N_{+,T})\sim C T^\frac{n(n+1)}{2}\;\;\hbox{as}\;\; T\rightarrow\infty
\end{equation}
for some $C>0$.

\subsection{Uniform distribution}

\hspace{1cm}\vspace{0.2cm}

\begin{proof}[Proof of Theorem \ref{th_frames0}]

We can write $v^0=g_0^{-1}e^0$ for some $g_0\in G$, where $e^0=(e_1,\ldots,e_n)$ is the standard frame.
Without loss of generality, $g_0^{-1}=k_0b_0$ for some $k_0\in K$, and $b_0\in B$  (see (\ref{eq_spiw})).

Let $Y=KBg_0$. By (\ref{eq_spiw2}), the product map $Y\times N_+^{g_0}\rightarrow G$ is a diffeomorphism.
For $g\in G$, denote $y_g=k_gb_gg_0$, the $Y$-component of $g$. The map
$$
\alpha:Y\rightarrow {\mathcal{W}}_n\simeq G/N_+^{g_0}:y\mapsto yv^0
$$
is a diffeomorphism. Denote by $\nu_1$ the measure on $Y$ which is the pull-back of the measure $\nu$ by the map $\alpha$.

For $g\in G$, $gv^0=y_gv^0$. This shows that $gv^0\in \Omega$ iff
$y_g\in \alpha^{-1} (\Omega)$.
Write
$$
n_g=\left(\begin{tabular}{cc}
$E$ & $l_g$\\
$0$ & $E$
\end{tabular}
\right).
$$
Then
$$
\|g\|=\|k_gb_gn_gg_0\|=\|b_gn_gb_0^{-1}\|=\left\|
\left(\begin{tabular}{cc}
$*$ & $b_g'l_g (^t b_0')$\\
$0$ & $*$
\end{tabular}
\right) \right\|.
$$
Thus,
$$
\|g\|^2=\|\hat c_g(n_gd_g)\|^2+e_g,
$$
where 
$$
c_g=b_g,\quad d_g=\left(\begin{tabular}{cc}
$(b_0')^{-1}b'_g$ & $0$\\
$0$ & $(^tb_0')(^t b_g')^{-1}$
\end{tabular}
\right),
$$
and $e_g$ is a continuous function depending only on the $B$-components of $g$.
We can use Proposition \ref{pro_assym} with $H=N_+$, $h_g=n_g$, and $m_g=1$.
Since $\Gamma\cdot v^0$ is dense in $\mathcal{W}_n$, $\Gamma N_+^{g_0}$ is dense in $G$. 
By (\ref{eq_last}), the condition (\ref{eq_H1}) holds for $H=N_+$.
Since $N_+$ is unipotent, the condition (\ref{eq_H2}) for $H=N_+$ holds too \cite{sh94}.

Applying Proposition \ref{pro_assym}, we get
\begin{eqnarray*}
N_T(\Omega,v^0)\sim\left(\frac{1}{\bar\mu(\Gamma\backslash G)}\int_{\alpha^{-1}(\Omega)}
\frac{d\nu_1(y)}{\Delta_{H}(c_y)}\right)\varrho(N_{+,T})
\end{eqnarray*}
as $T\rightarrow\infty$, where $\Delta_H$ is defined in (\ref{eq_DH}).
Thus, by (\ref{eq_last}),
$$
N_T(\Omega,v^0)\sim\lambda_{v^0}(\Omega)T^\frac{n(n+1)}{2}\;\;\hbox{as}\;\;T\rightarrow\infty,
$$
where
$$
\lambda_{v^0}(\Omega)=\frac{C}{\bar\mu(\Gamma\backslash G)}\int_{\Omega}
\frac{d\nu(x)}{\Delta_{N_+}(c_x)}.
$$
\end{proof}

\section*{Appendix}

\begin{proof}[Proof of Lemma \ref{l_frame_measure} (part 2)]
To find the constant $d_{n,l}$, we calculate measures of the set
$$
D=\{v=(v_1,\ldots, v_{l})\in {\mathcal{V}}_{n,l}: \|v_i\|<1,1\le i\le l\}.
$$
Denote by $V_k$ the Lebesgue measure of a $k$-dimensional unit  ball.
Recall that
\begin{equation} \label{eq_vball}
V_k=\frac{\pi^{k/2}}{\Gamma(1+k/2)}.
\end{equation}
Clearly,
\begin{equation} \label{eq_vol1}
\int_D dv=V_n^{l}=\frac{\pi^{nl/2}}{\Gamma(1+n/2)^{l}}.
\end{equation}
For $k\in K$, $n\in N_{l-}$, and $a\in A^o_{l-}$, $knae^0\in D$ iff $\|knae_i\|<1$
for $i=1,\ldots, l$. We have
\begin{equation}\label{eq_nnn}
\|kn^{_-}(t)a(s^{_-})e_i\|^2=\exp(2s_i^{_-})\left(1+t_{1i}^2+\cdots +t_{i-1 i}^2\right).
\end{equation}
Let as introduce new coordinates on $A^o_{l-}$: $a_i=\exp(s^{_-}_i)$, $1\le i\le l$.
The Haar measure on $A^o_{l-}$ (\ref{eq_da}) is given by $da=\prod_{i=1}^l\frac{da_i}{a_i}$.
By (\ref{eq_nnn}), the set of $(k,n^{_-}(t),a)\in K\times N_{l-}\times A^o_{l-}$ such that $kn^{_-}(t)ae^0\in D$
is described by conditions:
\begin{eqnarray*}
0<a_i<1\quad &&i=1,\ldots, l,\\
\|t_{*i}\|<\left(\frac{1-a_i^2}{a_i^2}\right)^{1/2}\quad &&i=2,\ldots, l.
\end{eqnarray*}
Thus,
\begin{eqnarray}
&&\nonumber\mathop{\int}_{kn^{_-}a^{_-}e^0\in D}\delta^-_l(a^{_-})^ndkdn^{_-}da^{_-}\\
\nonumber &=&\int_0^1\cdots\int_0^1\prod_{i=2}^{l}V_{i-1}\left(\frac{1-a_i^2}{a_i^2}\right)^{\frac{i-1}{2}}\delta_l^-(a^{_-})^n da^{_-}\\
&=&\nonumber\prod_{i=1}^{l} V_{i-1}\int_0^1 (1-a_i^2)^\frac{i-1}{2}a_i^{n-i}da_i\\
\nonumber &=&\prod_{i=1}^{l} \frac{V_{i-1}}{2}\int_0^1 (1-b_i)^\frac{i-1}{2} b_i^\frac{n-i-1}{2} db_i\\
&=&2^{-l}\prod_{i=1}^{l}V_{i-1}B\left(\frac{i+1}{2},\frac{n-i+1}{2}\right) \nonumber\\
&=&\frac{\pi^{l(l-1)/4}}{2^{l}\Gamma(1+n/2)^{l}}\prod_{i=1}^{l} \Gamma\left(\frac{n-i+1}{2}\right).\label{eq_vol2}
\end{eqnarray}
In the last step, we have used (\ref{eq_vball}) and the well-known identity for $\Gamma$-function and $B$-function.
Finally, by (\ref{eq_vol1}) and (\ref{eq_vol2}),
\begin{equation} \label{eq_c0}
d_{n,l}=\frac{2^{l}\pi^{l(2n-l+1)/4}}{\Gamma(n/2)\Gamma((n-1)/2)\cdots \Gamma((n-l+1)/2)}.
\end{equation}
\end{proof}

\begin{proof}[Proof of Lemma \ref{lem_BTasy}]
Let 
$$
d(\bar a)=\hbox{diag}\left(1,\ldots,1,a_{l+1},\ldots,a_n\right)
$$
for $\bar a=(a_{l+1},\ldots,a_n)\in\mathbb{R}^{n-l}$, and 
$$
\alpha(\lambda)=\hbox{diag}\left(1,\ldots,1,\lambda^{\frac{1}{n-l}},\ldots,\lambda^{\frac{1}{n-l}}\right).
$$
For $b=n^{_+}(t)d(\bar a)$, define
\begin{equation}\label{eq_La}
\Lambda (b)=\sum_{l<i\le n} a_i^2+\sum_{\max (i,l)<j} a_j^2t_{ij}^2.
\end{equation}
Note that $\Lambda(b)=\|b\|^2-l$.
Thus, it is enough to compute asymptotics of the function
$$
\phi(x)\stackrel{def}{=}\int_{\Lambda(b)<x} d\varrho_l(b).
$$ 
as $x\rightarrow\infty$. By Tauberian theorem \cite[Ch.~V, Theorem~4.3]{wid}, it can be deduced from asymptotics of
the function
\begin{equation} \label{eq_psi}
\psi(\lambda)\stackrel{def}{=}\int_0^\infty e^{-\lambda x}d\phi(x)=\int_{B^o_l} \exp\{-\lambda \Lambda(b)\} d\varrho_l(b).
\end{equation}
as $\lambda\rightarrow 0+$. It is more convenient to work with the function
\begin{equation} \label{eq_psi2}
\tilde{\psi}(\lambda)=\psi\left(\lambda^{\frac{2}{n-l}}\right).
\end{equation}
Let $\bar B^o_l=N_{l+}d({\mathbb{R}_+^d})=B^o_l\alpha(\mathbb{R}_+)$.
One can check that 
\begin{equation}\label{eq_new}
\int_{N_{l+}\times {\mathbb{R}_+^{n-l}}} f(nd(\bar a))dn^{_+}\frac{da_{l+1}}{a_{l+1}}\cdots\frac{da_{n}}{a_{n}} =\int_{B^o_l\times\mathbb{R}_+} f(b\alpha(\lambda))d\varrho_l(b)\frac{d\lambda}{\lambda}
\end{equation}
for $f\in L^1(\bar B^o_l)$. (In fact, each of the integral defines a right Haar measure on $\bar B^o_l$.)

Consider Mellin transform of the function $\tilde{\psi}$:
\begin{eqnarray*}
F(z)&=&\int_0^\infty \tilde\psi(\lambda)\lambda^z\frac{d\lambda}{\lambda}
\stackrel{(\ref{eq_psi})}{=}\int_{B^o_l\times\mathbb{R}_+}\exp\left\{-\Lambda(b\alpha(\lambda))\right\} \lambda^z d\varrho_l(b)\frac{d\lambda}{\lambda}\\
&\stackrel{(\ref{eq_new})}{=}&\mathop{\int}_{N_{l+}\times\mathbb{R}^{n-l}_+}\exp\left\{-\Lambda(n^{_+}d(\bar a))\right\} \left(\prod_{i=l+1}^n a_{i}\right)^z dn^{_+}\frac{da_{l+1}}{a_{l+1}}\cdots\frac{da_{n}}{a_{n}}\\
&\stackrel{(\ref{eq_La})}{=}&\mathop{\int}_{N_{l+}\times\mathbb{R}^{n-l}_+}\exp\left\{-\sum_{\max (i,l)<j}(a_jt_{ij})^2-\sum_{i=l+1}^na_i^2\right\}\\
&\times&\left(\prod_{i=l+1}^n a_{i}\right)^{z-1}\prod_{\max (i,l)<j}dt_{ij}\prod_{i=l+1}^n da_{i}.
\end{eqnarray*}
Using that $\int_\mathbb{R} e^{-x^2}dx=\sqrt{\pi}$, we get
\begin{eqnarray*}
F(z)=\pi^\frac{(n+l-1)(n-l)}{4}\int_{\mathbb{R}^{n-l}_+}\exp\left\{-\sum_{i=l+1}^n a_i^2\right\}
\left(\prod_{i=l+1}^n a_i^{z-i}\right) da_{l+1}\ldots da_n.
\end{eqnarray*}
Making substitution $u_i=a_i^2$, we get
\begin{equation}\label{eq_Fz}
F(z)=\frac{\pi^\frac{(n+l-1)(n-l)}{4}}{2^{n-l}}\prod_{j=l+1}^n \Gamma\left(\frac{z-j+1}{2}\right).
\end{equation}
By Mellin inversion formula, for sufficiently large $u$,
\begin{equation}\label{eq_fff}
\tilde\psi(\lambda)=\frac{1}{2\pi i}\int_{\hbox{\small Re}(z)=u} F(z)\lambda^{-z}dz.
\end{equation}
Since $\Gamma$-function decays fast on vertical strips, we can shift the line of integration
to the left. By (\ref{eq_Fz}), the first pole of $F(z)$ occurs at $z=n-1$. Therefore,
it follows from (\ref{eq_fff}) that
$$
\tilde\psi(\lambda)\sim \frac{\pi^\frac{(n+l-1)(n-l)}{4}}{2^{n-l-1}}\prod_{j=l+1}^{n-1}
\Gamma\left(\frac{n-j}{2}\right)\lambda^{-(n-1)}\;\;\hbox{as}\;\;\lambda\rightarrow 0+.
$$
By (\ref{eq_psi2}),
$$
\psi(\lambda)\sim \frac{\pi^\frac{(n+l-1)(n-l)}{4}}{2^{n-l-1}}\prod_{j=l+1}^{n-1}
\Gamma\left(\frac{n-j}{2}\right)\lambda^{-\frac{(n-1)(n-l)}{2}}\;\;\hbox{as}\;\;\lambda\rightarrow 0+.
$$
Finally, the asymptotic estimate for $\phi(x)$ as $x\rightarrow\infty$ follows from
Tauberian Theorem \cite[Ch.~V, Theorem~4.3]{wid}. We have
\begin{equation}\label{eq_gamnl}
\gamma_{n,l}=\frac{\pi^{(n+l-1)(n-l)/4}}{2^{n-l-1}\Gamma\left(\frac{(n-1)(n-l)}{2}+1\right)}
\prod_{j=l+1}^{n-1} \Gamma\left(\frac{n-j}{2}\right).
\end{equation}
This proves the lemma.
\end{proof}

The following lemma is used in the proof of Lemma \ref{lem_BTC2}.

\begin{lem}\label{lem_BTC}
For $C\in\mathbb{R}$ and $T>0$,
\begin{equation}\label{eq_rBTC}
\varrho_l(B^C_{l,T})=c_{n,l}\int_{A^C_T} \Big(T^2-N(s^{_+})-l\Big)^{\frac{(n-l)(n+l-1)}{4}}\hbox{\rm exp}\left\{\sum_{k=l+1}^n
(n-k)s^{_+}_k\right\}ds^{_+},
\end{equation}
where
\begin{eqnarray}
c_{n,l}&=&\frac{\pi^{(n-l)(n+l-1)/4}}{\Gamma(1+(n-l)(n+l-1)/4)},\nonumber\\
N(s^{_+})&=&\sum_{i=l+1}^n \exp\{2s^{_+}_i\}.\label{eq_Ns}
\end{eqnarray}
\end{lem}

\begin{proof}
Note that
$$
B^C_{l,T}=\left\{a(s^{_+})n(t): l+N(s^{_+})+\sum_{l<i<j\le n} \exp\{2s^{_+}_i\}t_{ij}^2+\sum_{1\le i\le l; l<j\le n} t_{ij}^2<T^2\right\}.
$$
We use the formula (\ref{eq_rho_l}) for $\varrho_l$ and make the change of variables 
$$t_{ij}\rightarrow \exp\{-s^{_+}_i\}t_{ij}$$
for $l<i<j\le n$.
The formula (\ref{eq_rBTC}) follows from the fact
that the volume of a unit ball in $\mathbb{R}^m$ is $\pi^{m/2}/\Gamma(1+m/2)$.
\end{proof}

\begin{proof}[Proof of Lemma \ref{lem_BTC2}]
For $i_0=l+1,\ldots, n-1$, put 
$$
A_{l+,T}^{i_0}=\{a(s^{_+})\in A_{l+,T}^o: s^{_+}_{i_0}\le C\}\;\;\hbox{and}\;\;B_{l,T}^{i_0}=\{a(s^{_+})n(t)\in B_{l,T}^o: s^{_+}_{i_0}\le C\}.
$$
We claim that $\varrho_l (B_{l,T}^{i_0})=o(\varrho (B_{l,T}^o))$ as $T\rightarrow\infty$.
If $a(s^{_+})\in A_{l+,T}^o$, then
$s^{_+}_i<\log T$ for every $i=l+1,\ldots, n$.
Then as in Lemma \ref{lem_BTC}, 
\begin{eqnarray*} 
\varrho_l(B_{l,T}^{i_0})&\le& c_{n,l}T^{\frac{(n-l)(n+l-1)}{2}}\int_{A_{l+,T}^{i_0}} \exp\left\{\sum_{k=l+1}^n
(n-k)s^{_+}_k\right\}ds^{_+} \\
&\ll& T^{\frac{(n-l)(n+l-1)}{2}} \prod_{{l<k<n},{k\ne i_0}} \int_{-\infty}^{\log T} \exp\{(n-k)s^{_+}_k\}ds^{_+}_k\\
&\ll& T^{(n-1)(n-l)-(n-i_0)}.
\end{eqnarray*}
Now the claim follows from (\ref{eq_BTasy}).
Since 
$$
B_{l,T}^o-B_{l,T}^C=\bigcup_{l<i_0<n} B_{l,T}^{i_0},
$$
we have
$$
\varrho_l (B_{l,T}^o-B_{l,T}^C)=o(\varrho_l (B_{l,T}^o))\;\;\hbox{as}\;\; T\rightarrow\infty.
$$
Therefore, $\varrho_l(B_{l,T}^C)\sim\varrho_l(B_{l,T}^o)$ as $T\rightarrow\infty$.
\end{proof}

\begin{proof}[Proof of Corollary \ref{th_frames}]
By Theorem \ref{th_dr}, $\Gamma v^0$ is dense in $\mathcal{V}_{n,l}$.
By Theorem \ref{th_frames00}, (\ref{eq_NTasy}) holds.
The volume of $\Gamma\backslash G$ was computed by Minkowski.
For the measure $\bar\mu$, we have 
$$
\bar \mu(\Gamma\backslash G)=2^{-(n-1)}\prod_{i=2}^n\pi^{-i/2}\Gamma(i/2)\zeta(i)
$$
(see \cite[Theorem 5.6]{sh00}).

Therefore,
$$
N_T(\Omega,v^0)\sim b_{n,l}\hbox{\rm Vol}(v^0)^{1-n}\left(\int_{\Omega} \frac{dv}{\hbox{\rm Vol}(v)}\right) T^{(n-1)(n-l)}\quad\hbox{as}\quad T\rightarrow\infty,
$$
where
\begin{equation}\label{eq_bnl}
b_{n,l}=\frac{a_{n,l}}{\bar\mu(\Gamma\backslash G)}=\frac{\pi^{n(n-l)/2}}{\Gamma\left(\frac{(n-1)(n-l)}{2}+1\right)\Gamma\left( \frac{n-l}{2}\right)}\prod_{i=2}^n\zeta(i)^{-1}.
\end{equation}
Here we used that $\Gamma(\frac{1}{2})=\sqrt{\pi}$.
\end{proof}

\end{document}